\documentclass[12pt]{article}
\usepackage[mathscr]{eucal}
\usepackage{amsbsy}
\usepackage{amssymb}
\usepackage{amsmath}
\usepackage{amsthm}
\usepackage{graphics}
\usepackage{epsfig}

\usepackage{latexsym}

\def\phi{\varphi}

\def\bbr{{\mathbb R}}

\def\bz{{\bf Z}}

\def\bz{{\mathbb Z}} 

\def\Bl1{{\bf 1}}
\def\B2{{\bf 2}}
\def\B0{{\bf 0}}

\def\a{\alpha}
\def\b{\beta}
\def\d{\delta}

\def\e{\varepsilon}
\def\g{\gamma}

\def\l{\lambda}

\def\s{\sigma}

\def\cd{{\mathscr D}}

\def\H{{\mathscr H}}
\def\=A8{\"o}

\def\wc{\stackrel{\cd}{\longrightarrow}}
\def\pc{\stackrel{P}{\longrightarrow}}

\def\fdd{\stackrel{f.d.d.}{\longrightarrow}}

\newcommand{\dd}{{\rm d}}
\newcommand{\sign}{{\rm sign} }
\newcommand{\ri}{{\rm i} }
\newcommand{\PP}{\mathbf{P}}

\newcommand{\beq}{\begin{equation}}
\newcommand{\eeq}{\end{equation}}
\newcommand\beqn{\begin{displaymath}}  
\newcommand\eeqn{\end{displaymath}}

  \newcommand{\RR}{\mathbb{R} }
\newcommand{\abs}[1]{\left|#1\right|}

\newcommand{\du}{{\rm d}u }

\newcommand{\ind}[1]{\mathbbm{1}_{#1}}
\newcommand{\sv}[1]{\left\lfloor #1 \right\rfloor}
\newcommand{\halmos}{\vspace{3mm} \hfill \mbox{$\Box$}\\[2mm]}
\theoremstyle{plain}
\newtheorem{teo}{Theorem}

\newtheorem{prop}[teo]{Proposition}
\theoremstyle{definition}
\newtheorem{definition}[teo]{Definition}

\newtheorem{remark}[teo]{Remark}

\usepackage{bbm}
\begin{document}

\title{Functional limit theorems for linear processes with tapered innovations and filters
\footnotemark[0]\footnotetext[0]{ \textit{Short title:}
linear processes with tapered innovations }
\footnotemark[0]\footnotetext[0]{%
\textit{MSC 2010 subject classifications}. Primary 60G99, secondary
60G22, 60F17 .} \footnotemark[0]
\footnotemark[0]\footnotetext[0]{ \textit{Corresponding author:}
Vygantas Paulauskas, Department of Mathematics and
 Informatics, Vilnius university, Naugarduko 24, Vilnius 03225, Lithuania,
 e-mail:vygantas.paulauskas@mif.vu.lt}
}


\author{ Vygantas Paulauskas$^{\text{\small 1}}$ \\
{\small $^{\text{1}}$ Vilnius University, Department of Mathematics
and
 Informatics,}\\}



\maketitle

\begin{abstract}

In the paper we consider the partial sum process $\sum_{k=1}^{[nt]}X_k^{(n)}$, where $\{X_k^{(n)}=\sum_{j=0}^{\infty} a_{j}^{(n)}\xi_{k-j}(b(n)), \ k\in \bz\},\ n\ge 1,$ is a series of linear processes with tapered filter $a_{j}^{(n)}=a_{j}\ind{[0\le j\le \l(n)]}$ and  heavy-tailed tapered innovations $\xi_{j}(b(n), \ j\in \bz$. Both tapering parameters  $b(n)$ and $\l(n)$ grow to  $\infty$ as $n\to \infty$. The limit behavior of the partial sum process (in the sense of convergence of finite dimensional distributions) depends on the growth of these two tapering parameters and dependence properties of a linear process with non-tapered filter $a_i, \ i\ge 0$ and non-tapered innovations. We consider the cases where $b(n)$ grows relatively slow (soft tapering) and rapidly (hard tapering), and all three cases of growth of $\l(n)$ (strong, weak, and moderate tapering).

\end{abstract}
\section{Introduction }

Linear processes are widely used both in theoretical and applied probability. Let us define
\begin{equation}\label{linpr1}
X=\{X_k, \ k\in \bz \}, \ X_k=\sum_{j=0}^\infty a_j\xi_{k-j}.
\end{equation}
Here the filter $\{a_j, j\ge 0\}$ and innovations $\xi_i, i\in \bz$, independent and identically distributed (i.i.d.) random variables, are such that the linear random process in (\ref{linpr1}) is well defined.  By means of linear processes, depending on moments of innovations, we can model stationary sequences with finite or infinite variance, while the properties of a filter of a linear process allow to model different dependence of a sequence $\{X_k\}$, namely, long-range, short-range, and negative dependence (we shall use the abbreviations LRD, SRD, and ND, respectively). Let us denote
\begin{equation}\label{sum}
S_{n}=S_{n}(X):=\sum_{k=1}^{n} X_{k} \ {\rm and} \ S_{n}(t; X)=S_{[nt]}(X), \quad t\ge 0.
\end{equation}
The asymptotic behavior of partial sum processes $S_{n}(t; X)$ is well investigated,  starting with pioneering works  of Davydov \cite{dav1970} (the case of finite variance) and Astrauskas \cite{Astrauskas} (the case of infinite variance) and with a big list of subsequent papers, devoted to limit theorems for $S_{n}(t; X)$.

In many fields, especially connected with applications, we face random quantities distributed according power-law. For example, it is well-known that many  natural hazards,  such as earthquakes, rock falls, landslides, riverine floods, tsunami, wildfire exhibit power-law behavior, see, for example, \cite{Geist}. At the same time we must admit that essentially real world is bounded, i.e., the quantities, which we are interested in, have some bounds or at least have much lighter tails. Good example is seismology. For many years the scalar value of seismic moment of an earthquake was modeled by Pareto distribution with exponent $\a$. Theoretical reasoning, based on the theory of branching processes, predicts that this exponent is universal (for many types of earthquakes - ground, oceanic, shallow, etc.) constant equal $1/2$, although estimations from real data usually give a little bit bigger value, see discussion in \cite{Kagan}. Another important fact was that empirical data of earthquakes showed that while seismic moment follows power-law in quite big range of values, the largest values of samples demonstrate much lighter tails, see Figure 1 in \cite{Kagan}. Therefore it was suggested to apply the exponential taper to Pareto distribution. There are  several ways to define tapering of a random variable, for example exponentially tapered Pareto distribution was introduced by V. Pareto himself, see, for example, \cite{Vaiciulis}. We shall use definition of tapering which was used in \cite{Chakrabarty}.  If $\xi$ has a heavy-tailed distribution with the tail index $\a$ and $R$ is a non-negative random variable, independent from $\xi$ and having the light tail, then tapered random variable $\xi(b)$ with the taper $R$ and tapering parameter $b>0$ is defined as follows
$$
\xi(b)=\xi \ind{[|\xi|<b]}+\frac{\xi}{|\xi|}(b+R)\ind{[|\xi|\ge b]}.
$$
Usually $R$ is exponential law, but it can be degenerated at zero random variable, then $\xi (b)$ will be truncated (at level $b>0$) random variable. An example, illustrating the last possibility can  be provided from computer science. In many problems it is assumed that the file size in networks has heavy-tailed distribution, on the other hand there is so-called the File Allocation Table (FAT, see Microsoft Knowledge Base Article 154997 (2007) or Wikipedia), used on most computer systems, which allows the
largest file size to be 4GB (more precisely, $2^{32}-1$ bytes). There are more examples where tapered or truncated  heavy-tailed distributions are used to model various processes, we can recommend papers \cite{Aban} and \cite{Meerschaert} containing big lists of references with such examples.
Tapering parameter can be dependent on $n$, and if $b_n$ is growing to infinity comparatively slowly, we say that we have hard tapering, if it grows rapidly, then we have soft tapering and between them there is intermediate tapering  (for strict definitions see Definition \ref{def1} in Section 3.1). The next natural step is to consider the family of linear processes $X^{(n)}=\{X_k^{(n)}\}$, indexed by $n$ and for each $n$ taking as innovations a sequence of i.i.d.  copies of a tapered heavy-tailed random variable $\xi(b_n)$ and to investigate  $S_{n}(t; X^{(n)})$. In the case where for each $n$ random variables $X_k^{(n)}, \ k\in \bz $ are i.i.d. (this will be if $a_0=1, a_j=0$ for $j\ge 1$) such answer for $S_{n}(1; X^{(n)})$ (even in multidimensional case) is given in \cite{Chakrabarty}. Namely, in the case of hard tapering we have Gaussian limit law for appropriately normalized $S_{n}(1; X^{(n)})$, while in the case of soft tapering limit for $S_{n}(1; X^{(n)})$ will be the same as for sum $\sum_{k=1}^{n} \xi_{k}$ with the same centering and normalizing sequences. In the case  of intermediate tapering in \cite{Chakrabarty} it was shown that limit law is some infinitely divisible law.
In papers \cite{Paul21}, \cite{Paul23}, \cite{Paul22}, and \cite{PaulDam5} limit theorems for linear random processes and linear random fields (with a specific structure of a filter) with tapered innovations were considered. In the case of linear random processes with tapered innovations we got   almost the final result, see Theorem 3 in \cite{PaulDam5}. We use the word "almost", since it remains  as open problem  what is a limit process in the case of intermediate tapering.

Linear processes are formed by the means of a filter and innovations. In \cite{Sabzikar} linear processes with tapered filters were introduced and limit theorems  for sums of values of such processes
were investigated.
In \cite{Paul21} it was mentioned that it is possible to consider not only tapered innovations, but together  tapered innovations and a filter, and in this paper we shall investigate this possibility. Now we introduce the notion of a linear process with a tapered filter, and our definition will be a little bit more general than definition in \cite{Sabzikar}. By a  taper for the filter of a linear process (\ref{linpr1}) we call a sequence $f(j)>0, \ j\ge 0,$ which is rapidly
decreasing as $j\to \infty$ and such that $\sum_{j=0}^\infty |a_jf(j)|<\infty$.
  Then a linear process with a tapered filter is defined as $X_k=\sum_{j=0}^\infty a_jf(j)\xi_{k-j}, \ k\in \bz.$
Such linear process with a tapered filter, independently from initial filter $\{a_j\}$, is always with short-range dependence, due to the condition  $\sum_{j=0}^\infty |a_jf(j)|<\infty$, therefore, considering sums of values of a linear process, it is interesting to consider family of linear processes with tapers $f_n(j)$, depending on $n$. Additionally, it is natural to require that $f_n(j)$ is close to $1$ for $j=0, 1, \dots, k_n$ with some  $k_n \to \infty$, as
 $n\to \infty$.
Thus, we shall consider a family of linear processes
\begin{equation}\label{trunkfil}
X^{(n)}=\{X_k^{(n)} \}, \quad {\rm where} \ \  X_k^{(n)}=\sum_{j=0}^\infty a_jf_n(j)\xi_{k-j}, \ k\in \bz.
\end{equation}
In \cite{Sabzikar} the following family of tapers $f_n(j)=\exp (-j/\l(n))$ with $\l(n)\to \infty$, as $n\to \infty$, was used. If there exists the limit
$$
\lim_{n\to \infty}\frac{n}{\l(n)}=\l_*,
$$
 then the family $X^{(n)}$ is  called {\it strongly, weakly, or moderately} tapered, if $\l_*=\infty, \l_*=0,$ or $0<\l_*<\infty,$ respectively. Let us note that in \cite{Sabzikar} the term "tempered"was used instead of "tapered" and the taper was of the form $\exp (-j{\bar \l}(n))$ with ${\bar \l}(n)\to 0$ as $n\to \infty$. Since we shall deal with tapered innovations and filter together, for us it is more convenient to have both parameters of tapering unboundedly increasing. Also we  take a  different family of tapers, namely, $f_n(j)=\ind{[0\le j\le \l(n)]},$ and for initial filter $\{a_j\}$ the  tapered filter will be denoted
\begin{equation}\label{newtaper}
{\bar a}_i= {\bar a}_i^{(n)}=\left \{\begin{array}{ll}
              a_i, & {\rm if} \ 0\le i\le \l(n), \\
              0, & {\rm if} \ i> \l(n).
               \end{array} \right.
\end{equation}
We think that this new taper (\ref{newtaper}) (which essentially is the truncation) is more natural, than that used in \cite{Sabzikar}. If tapering by exponential taper of innovations, which are random variables,  is quite natural (in many cases more natural than truncation), for filter coefficients (non random quantities) more natural is truncation, since in real life we always have only finite number of filter coefficients, only this finite number can be big. This tapering gives one more meaning for our family of random variables (\ref{trunkfil}). Since $X_k^{(n)}=\sum_{j=0}^m a_j^{(n)}\xi_{k-j}, \ k\in \bz$ with $m=\l(n)$, thus we have series of $m$-dependent random variables $X_k^{(n)}, \ k\in \bz$ with $m \to \infty$. Sums of $m$-dependent random variables with unboundedly increasing $m$ appears in many problems of theoretical and applied probability. It should  be interesting to compare these results (for example, results from \cite{Romano}) with results of our paper.
%
 We shall take $\l(n)=cn^{\g_1} \quad {\rm with}\ \ 0<c<\infty, \ 0<\g_1<\infty.$
We have strongly, weakly, or moderately tapered family of linear processes (\ref{trunkfil}) (with taper from (\ref{newtaper})), if $0<\g_1<1, \ \g_1>1,$ or $\g_1=1$, respectively. It is clear that in the cases of strong and weak tapering the constant $c$ does not play any role, therefore in the cases $\g_1\ne 1$ we shall take $c=1$. In the case of moderate tapering  $\l_*=c^{-1}$, and this constant is important.

The rest of the paper contains two sections. In Section 2 we consider the case where only the filter is tapered, in Subsection 2.1 the Gaussian case is considered and in Subsection 2.2 we investigate stable case.
 Results, obtained in these subsections  are similar to those obtained in Section 4 in \cite{Sabzikar}. The only difference is in the case of moderate tapering, limit processes in this case in Theorems \ref{thm1} and \ref{thm3} is different from those in \cite{Sabzikar} and gives new examples of Gaussian and stable processes. This is discussed in Subsection 2.3. Here it is appropriate to mention that there is a small problem in \cite{Sabzikar}: the statement of Proposition 4.2, which is given without proof and which is the main step in proving Theorem 4.3 (limit theorem for linear processes with exponentially tapered filter), in the case $1<\a<2$ is incorrect, D. Surgailis himself constructed counterexample (personal communication).

Finally, in Section 3, divided into two subsections, we consider family of linear processes with tapered innovations and filter
\begin{equation}\label{trunkfilinn}
{\bar X}^{(n)}=\{{\bar X}_k^{(n)} \}, \quad   {\bar X}_k^{(n)}=\sum_{i=0}^\infty {\bar a}_i^{(n)}\xi_{k-i}(b(n)), \ k\in \bz,
\end{equation}
where the filter  is from (\ref{newtaper}) and $\xi_{k}(b(n)), \ k\in \bz,$ are tapered innovations, defined in Section 3, see (\ref{paretodf}) and (\ref{tapparetodf}). We  take, as in \cite{Paul21},  $b(n)=n^\g$ and we  have hard, soft or intermediate tapering, if $0<\g<1/\a, \ \g>1/\a$, or $\g=1/\a$, respectively. We left the same terminology, used in \cite{Sabzikar} for tapering of the filter, since now we can say, for example, that we consider a linear process with weak and hard tapering, understanding that we have weak tapering of the filter and  hard tapering of innovations.

\section{Limit theorems for linear processes with tapered filters}\label{sec1}

\subsection{The Gaussian case}\label{subsec1}

In this section we consider the family of linear processes
$$
X^{(n)}=\{X_k^{(n)} \}, \quad   X_k^{(n)}=\sum_{j=0}^\infty {\bar a}_j^{(n)}\e_{k-j}, \ k\in \bz,
$$
where ${\bar a}_i^{(n)}$ is defined in (\ref{newtaper})
 and $\e_{k}, \ k\in \bz,$ are i.i.d. random variables with $E\e_1=0, \ E\e_1^2=1.$
  Let us denote
$$
Z_n(t)=A_n^{-1} S_n(t, X^{(n)} ),
$$
where $A_n$ is a normalizing sequence for $ S_n(1, X^{(n)})$ and $S_n(t, X)$ is defined in (\ref{sum}).

 As in previous our papers \cite{Paul21}, \cite{Paul22}, and \cite{PaulDam5}, we assume that
\begin{equation}\label{cond1}
a_n\sim n^{-\b},
\end{equation}
where $\b >1/2$ (this condition ensures the correctness of definition of a linear process $X=\{X_k, \ k\in \bz\}$ with non-tapered filter; for the existence of $X_k^{(n)}$ for each $n$ this condition is irrelevant, later we shall discuss what is happening when $\b$ is even negative). Together with (\ref{cond1}) we shall consider three main sets of conditions, giving rise to three different types of dependence of the process $X$.

(i) $1/2<\b<1$ - the case of LRD;

(ii) $\b>1$ and $\sum_{i=0}^\infty a_i\ne 0$, the case of SRD;

(iii) $1<\b<3/2$ and $\sum_{i=1}^\infty a_i=0$, the case of ND.

If the linear process $X_k=\sum_{j=0}^\infty  a_j\e_{k-j}$ is with non-tapered filter, then it is well-known that
$\{n^{-H}S_n(t, X)\}\fdd \{B_H (t)\}$
where
$H=\frac{3}{2}-\b $, in  the  cases  (i)  and (iii), and $H=1/2$ in the case (ii).
Here $B_H$  is fractional Brownian motion (FBM) with Hurst parameter $0<H<1$, and $\fdd $ stands for the  convergence of finite-dimensional distributions (f.d.d.). For $H=1/2$ we  have Brownian motion (BM) $B(t):=B_{1/2}(t)$. We denote by $\wc$ the weak convergence of random processes in the Skorohod space $D[0,1]$ equipped with $J_1$-topology, see \cite{Billingsley}

Since for a filter we have three cases of tapering and three cases of dependence, thus we  consider nine cases. We shall number them by index $j=1,\dots , 9$ in the following way: $j=1, 2, 3$ we attribute to strong tapering ($0<\g_1<1$) and three dependence cases (i), (ii), and (iii), respectively, indices $j=4, 5, 6$ are attributed to week tapering, and $j=7, 8, 9$ are given to the case of moderate tapering and three dependence types. For example, $j=8$ means that we consider moderate tapering of  a filter and SRD. Therefore, we shall use the notations $Z_n^{(j)}(t), (A_n^{(j)})^2, S_n^{(j)}(t, X^{(n)})$, but in cases where it is clear which index $j$ is considered, we shall skip this index from the notation.

\begin{teo}\label{thm1} Suppose that   there exists $0<\d\le 1$ such that $E|\e_1|^{2+\d}<\infty$. Then, for all $j=1,\dots , 9$, we have
\begin{equation}\label{fddconv}
\{Z_n^{(j)}(t)\}\fdd \{U^{(j)}(t)\},
\end{equation}
where normalizing sequences $(A_n^{(j)})$ and Gaussian limit processes $U^{(j)}(t)$ are  defined by means of their covariance functions, which  are given in Proposition \ref{prop1}. Particularly, $U^{(j)}(t)=B(t)$ for $j=1, 2, 3,5, 8,$ and $U^{(j)}(t)=B_H(t)$, \ $H=3/2-\b$,  for   $j=4, 6$. Processes $U^{(7)}(t)$ and $U^{(9)}(t)$ will be discussed in Section 2.3.

Convergence in (\ref{fddconv}) can be strengthen to
\begin{equation}\label{wcconv}
\{Z_n^{(j)}(t)\}\wc \{U^{(j)}(t)\}
\end{equation}
without any additional condition for $1\le j \le 8, j\ne 6$ while for $j=6, 9$ we need stronger moment condition $p=2+\d>2(3-2\b)^{-1}$.
\end{teo}
\begin{remark}\label{rem1} Most probably (\ref{fddconv}) in Theorem \ref{thm1} holds for $\d=0$, if in the proof we should use Lindeberg type condition. But it will require rather complicated calculations in all nine cases. On the other  hand, in \cite{Paul21} and \cite{Paul22},
 we used Lyapunov fractions of order $2+\d$ and now we were able to use them. And the main reason to leave moment condition with $\d>0$ was that this condition gives us stronger convergence in almost all cases of tapering.
\end{remark}

\begin{remark}\label{rem2} We included the strengthening (\ref{wcconv}) only in Theorem \ref{thm1}, the same strengthening can be added in Theorem \ref{thm2}, since the proof of tightness of distributions of $Z_n$ and ${\bar Z}_n$ easily follows from Propositions \ref{prop1} and \ref{prop4}, respectively. This strengthening in the stable case (Theorems \ref{thm3} and \ref{thm4}) is more complicated and is postponed for subsequent paper.
\end{remark}

{\it Proof of Theorem \ref{thm1}}.
In order to prove (\ref{fddconv}) for $Z_n(t)$ we shall use the same scheme of the proof as in \cite{Paul21}: we calculate ${\rm Var}S_n(t, X^{(n)})$, then we get $A_n^2={\rm Var}S_n(1, X^{(n)})$,  find a limit $\lim_{n\to \infty}{\rm Var}Z_n(t)$, and, finally, prove that f.d.d. of $Z_n(t)$ are asymptotically Gaussian.

\begin{prop}\label{prop1} For all $j=1,\dots , 9$ and $s, t >0$ we have
\begin{equation*}\label{VarZnfinal}
\lim_{n\to \infty}{\rm Var}Z_n^{(j)}(t)= W^{(j)}(t),
\end{equation*}
\begin{equation*}\label{CovZnfinal}
\lim_{n\to \infty}{\rm Cov}\left (Z_n^{(j)}(t), Z_n^{(j)}(s)\right )= W^{(j)}(t)+W^{(j)}(s)-W^{(j)}(|t-s|),
\end{equation*}
where
$$
W^{(j)}(t)= t^{2H(j)}, \ {\rm for} \  j=1, 2, 3, 4, 5, 6, 8,
$$
$$
W^{(7)}(t)= t^{2H(7)}C_{10}(t, \b, c), \ \ W^{(9)}(t)= t^{2H(9)}C_{17}(t, \b, c).
$$
Here
$$
H(j)=1/2,\ {\rm for} \  j=1, 2, 3,  5,  8, H(j)=3/2-\b, \ {\rm for} \  j=4, 6, 7, 9,
$$
 and $C_{10}(t, \b, c)$ and  $C_{17}(t, \b, c)$ are defined in (\ref{VarZn3}) and (\ref{VarZn3iii}), respectively. Normalizing sequences are defined as follows: $(A_n^{(j)})^2=\left (\sum_{i=0}^\infty a_{i} \right )^2 n$ for $j=2, 5, 8,$ \ $(A_n^{(j)})^2= C(\b)n^{1+2\g_1(1-\b)}$ for $j=1, 3$, \ $(A_n^{(j)})^2= C(\b)n^{2H(j)}$ for $j=4, 6,$ and \ $(A_n^{(j)})^2= C(\b, c)n^{2H(j)}$ for $j=7, 9.$
\end{prop}

{\it Proof of Proposition \ref{prop1}}.
Using formulae (2.2) and (2.3) from \cite{Paul21} and recalling that $E\e_1^2=1$  we can write
$$
S_n(t, X^{(n)})=\sum_{i=-\infty}^{[nt]} d_{n,i,t}\e_{i}, \ \ {\rm Var}S_n(t, X^{(n)})=\sum_{i=-\infty}^{[nt]}| d_{n,i,t}|^2,
$$
where  $d_{n,i,t}=\sum_{k=1}^{[nt]}{\bar a}_{k-i}$ for $i\le 0$ and $d_{n,i,t}=\sum_{k=i}^{[nt]}{\bar a}_{k-i}$ for $i> 0$. As in \cite{Paul21} we can write
\begin{equation*}\label{sumdnj}
\sum_{i=-\infty}^{[nt]} |d_{n,i,t}|^2=V_1(t) +V_2(t):=\sum_{i=-\infty}^0 |d_{n,i,t}|^2 + \sum_{i=1}^{[nt]} |d_{n,i,t}|^2.
\end{equation*}
To find the asymptotic of $V_1(t), V_2(t)$ with respect to $n$ for all $j=1,\dots , 9$, as in \cite{Paul21} instead of condition (\ref{cond1}), which means that $a_n=n^{-\b}(1+\d(n))$, where $\d(n)\to 0$, as $n\to \infty$, we can simply assume $a_n=n^{-\b}$ for $n\ge 1$. Although in the proof of this step some changes are needed, since in \cite{Paul21} the non-tapered filter was considered, but since $\l(n)\to \infty$ in all cases of tapering, there is no principal difficulties. Also, since the detailed proof of Theorem \ref{thm1} is in our preprint \cite{Paul24}, we provide here the proof in the case $\g_1=1$ only, since only in this  case, namely, for $j=7, 9$, we get some new limit Gaussian processes. 

 We start with the case $j=8$ of moderate tapering. Let us denote $m=[nt], m_1=[cn]$.  Now  $m_1\le m$ or $m_1\ge m$ if $c<t$ or $c>t$, respectively. Let us consider the case $0<t<c$, i.e., $m\le m_1$. We have (skipping the index $j=8$)
$$
V_2(t)=\sum_{k=1}^m \left (\sum_{i=0}^{m-k} a_{i}\ind{[0\le i\le m_1]} \right )^2=\sum_{k=1}^m \left (\sum_{i=0}^{m-k} a_{i}\right )^2
$$
and
\begin{equation}\label{asV2A}
\frac{V_{2}(t)}{n} \to t\left (\sum_{j=0}^\infty a_j \right)^2.
\end{equation}
Denoting ${\tilde n}=m_1-m$, we have
$$
V_1(t)=\sum_{i=0}^\infty \left (\sum_{k=1}^{m} a_{k+i}\ind{[0\le k+i\le m_1]} \right )^2=V_{1,1}(t)+V_{1,2}(t),
$$
where
$$
V_{1,1}(t)=\sum_{i=0}^{{\tilde n}} \left (\sum_{k=1}^{m} a_{i+k} \right )^2, \quad V_{1,2}(t)=\sum_{i={\tilde n}+1}^{m_1-1} \left (\sum_{k=1}^{m_1-i} a_{i+k} \right )^2.
$$
Denoting $b_i=\left (\sum_{k=1}^{m} a_{k+i} \right )^2$, we see that $b_i\to 0$ as $i\to \infty$, therefore
taking into account the relation ${\tilde n}=n(c-t)(1+o(1))$, it is not difficult to infer that
\begin{equation}\label{asV11tA}
\frac{V_{1,1}(t)}{n} \to 0.
\end{equation}

Using (\ref{cond1}) and evident estimate $m_1/n <c+1$ we can estimate
$$
  \frac{1}{n} \left (\sum_{k=1}^{m_1-i} a_{i+k} \right )^2 \le \frac{1}{n}\sum_{k=1}^{m_1-i} a_{i+k} \sum_{i=1}^{\infty} |a_{i}|\le  C|a_{i+1}|,
$$
and  we get
\begin{equation}\label{asV12tA}
\frac{V_{1,2}(t)}{n}\le C \sum_{={\tilde n}+1}^{m_1-1}|a_{i+1}|\to 0, \quad {\rm as} \ n\to \infty.
\end{equation}
Collecting (\ref{asV2A})-(\ref{asV12tA}), we have
\begin{equation}\label{asV1+V2}
\frac{V_{1}(t)+V_2(t)}{n} \to t\left (\sum_{j=0}^\infty a_j \right)^2.
\end{equation}
In the case $0<c\le t$, i.e., $m \le m_1$, in a similar way (only now denoting ${\tilde n}=m-m_1$) we get (\ref{asV1+V2}). Therefore, in the case $j=8$, taking $\left (A_n^{(8)}\right )^2=n\left (\sum_{j=0}^\infty a_j \right)^2 $, we get the  relation
\begin{equation*}\label{VarZn}
\lim_{n\to \infty}{\rm Var}Z_n^{(8)}(t)=t.
\end{equation*}

Let us consider the case of moderate tapering and LRD, the case $j=7$. We start assuming that $0<t\le c$, i.e., $m\le m_1$. In a similar way as in the case  $j=8$, only taking into account that now $m_1/m \to c/t$,  as $n\to \infty$, we  can get
\begin{equation}\label{asV1mod}
\frac{V_{1}(t)}{n^{3-2\b}} \to t^{3-2\b}\left (C_{1}(t, \b, c)+ C_{2}(t, \b, c) \right ),
\end{equation}
where
$$
C_{1}(t, \b, c)=\int_0^{(c/t)-1}\left (\int_0^{1}(x+y)^{-\b} dx\right )^2 dy,
$$
$$
C_{2}(t, \b, c)=\int_{(c/t)-1}^{(c/t)}\left (\int_0^{(c/t)-y}(x+y)^{-\b} dx\right )^2 dy.
$$
Similarly we get
\begin{equation}\label{asV2mod}
\frac{V_{2}(t)}{n^{3-2\b}} \to t^{3-2\b}C_{3}(\b), \ \ {\rm where} \ \ C_{3}(\b)=\int_0^1\left (\int_0^{1-y}x^{-\b} dx\right )^2 dy.
\end{equation}
From (\ref{asV1mod}) and  (\ref{asV2mod}) we have
\begin{equation}\label{asV1+V2mod}
\frac{V_1(t)+V_{2}(t)}{n^{3-2\b}} \to t^{3-2\b}C_{4}(t, \b, c),
\end{equation}
where
$$
C_{4}(t, \b, c)=C_{1}(t, \b, c)+ C_{2}(t, \b, c) + C_{3}(\b).
$$
In the case $0<c<t$, i.e., $m_1\le m$, in a similar way we can get
\begin{equation}\label{asV1+V2modA}
\frac{V_1(t)+V_{2}(t)}{n^{3-2\b}} \to t^{3-2\b}C_{8}(t, \b, c),
\end{equation}
where
$$
C_{8}(t, \b, c)=C_{5}(t, \b, c)+C_{6}(t, \b, c)+C_{7}(t, \b, c),
$$
$$
C_{5}(t, \b, c)=\int_0^{1-(c/t)}\left (\int_0^{(c/t)}x^{-\b} dx\right )^2 dy, \ C_{6}(t, \b, c)=\int_{1-(c/t)}^{1}\left (\int_0^{1-y}x^{-\b} dx\right )^2 dy,
$$
$$
C_{7}(t, \b, c)=\int_0^{(c/t)}\left (\int_0^{(c/t)-y}(x+y)^{-\b} dx\right )^2 dy.
$$
Combining (\ref{asV1+V2mod}) and (\ref{asV1+V2modA}) we have
\begin{equation}\label{asV1+V2modB}
\frac{V_1(t)+V_{2}(t)}{n^{3-2\b}} \to t^{3-2\b}C_{9}(t, \b, c),
\end{equation}
where $C_{9}(t, \b, c)=C_{4}(t, \b, c)$, for $0<t\le c,$ and $C_{9}(t, \b, c)=C_{8}(t, \b, c)$, for $0<c<t.$
There is continuity of $C_{9}$ at point $t=c$, since $C_{1}(c, \b, c)=C_{5}(c, \b, c)=0$ and it is easy to verify that $ C_{2}(c, \b, c) + C_{3}(\b)=C_{6}(c, \b, c)+C_{7}(c, \b, c)$.
Now we can choose the normalizing sequence $A_n^2$, to this aim we put $t=1$ in (\ref{asV1+V2modB}). It is easy to see that
$C_{9}(1, \b, c)= C_{4}(1, \b, c)$ if $ 1\le c,$ and $C_{9}(1, \b, c)= C_{8}(1, \b, c)$ if $ o<c<1.$
Therefore, taking $\left (A_n^{(7)}\right )^2= C_{9}(1, \b, c)n^{3-2\b}$ we get
\begin{equation}\label{VarZn3}
\lim_{n\to \infty}{\rm Var}Z_n^{(7)}(t)=t^{3-2\b}C_{10}(t, \b, c), \ {\rm where} \ \  C_{10}(t, \b, c)=\frac{C_{9}(t, \b, c)}{C_{9}(1, \b, c)}.
\end{equation}

It remains the case $j=9$, the case of moderate tapering and ND. The case of  ND is similar to LRD, only now we use the equality $\sum_{i=0}^{k} a_{i}=-\sum_{i=k+1}^\infty a_{i}$, therefore we provide only the final result. In the case $0<t<c,  m\le m_1, $ we get
\begin{equation}\label{asV1+V2modiii}
\frac{V_1(t)+V_{2}(t)}{n^{3-2\b}} \to t^{3-2\b}C_{11}(t, \b, c),
\end{equation}
where $C_{11}(t, \b, c)=C_{1}(t, \b, c)+ C_{2}(t, \b, c) + C_{12}(\b).$
In the case $0<c<t,  m_1\le m, $ we get
\begin{equation}\label{asV1+V2modAiii}
\frac{V_1(t)+V_{2}(t)}{n^{3-2\b}} \to t^{3-2\b}C_{13}(t, \b, c),
\end{equation}
where
$C_{13}(t, \b, c)=C_{7}(t, \b, c)+C_{14}(t, \b, c)+C_{15}(t, \b, c).$
The expressions of these above written constants are as follows:
$$
 C_{12}(\b)=\int_0^1\left (\int_{1-y}^\infty x^{-\b} dx\right )^2 dy   , \  C_{14}(\b)=\int_{1-(c/t)}^{1}\left (\int_{1-y}^\infty x^{-\b} dx\right )^2 dy,
$$
$$
C_{15}(t, \b, c)=\int_0^{1-(c/t)}\left (\int_{(c/t)}^\infty x^{-\b} dx\right )^2 dy.
$$
Combining (\ref{asV1+V2modiii}) and (\ref{asV1+V2modAiii}) we have
$$
\frac{V_1(t)+V_{2}(t)}{n^{3-2\b}} \to t^{3-2\b}C_{16}(t, \b, c),
$$
where
$C_{16}(t, \b, c)=C_{11}(t, \b, c) $ if $ 0<t\le c,$ and $C_{16}(t, \b, c)=C_{13}(t, \b, c) $ if $ 0<c< t.$
Taking $\left (A_n^{(9)}\right )^2= C_{16}(1, \b, c)n^{3-2\b}$, we get
\begin{equation}\label{VarZn3iii}
\lim_{n\to \infty}{\rm Var}Z_n^{(9)}(t)=t^{3-2\b}C_{17}(t, \b, c), \ {\rm where} \ \  C_{17}(t, \b, c)=\frac{C_{16}(t, \b, c)}{C_{16}(1, \b, c)}.
\end{equation}
\halmos

The second step in the proof of (\ref{fddconv}) is to prove that f.d.d. of $Z_n^{(j)}$ are asymptotically normal.
We recall that $m=[nt]$.
\begin{prop}\label{prop2} If $E\e_1=0, E|\e_1|^{2}=1$ and there exists $0<\d\le 1$ such that $E|\e_1|^{2+\d}<\infty$, then, for all $j=1, 2, \dots, 9,$ as $n\to \infty$,
\begin{equation}\label{Liapunov}
L^{(j)}(2+\d,n, t):=\frac{\sum_{k=-\infty}^m |d_{n,k,t}^{(j)}|^{2+\d}E|\e_1|^{2+\d}}{\left (\sum_{k=-\infty}^n (d_{n,k,1}^{(j)})^2\right )^{(2+\d)/2}} \to 0.
\end{equation}
\end{prop}
{\it Proof of Proposition \ref{prop2}}. From the proof of Proposition \ref{prop1}, for any fixed $t$ and for all $j=1, 2, \dots, 9,$ we can estimate
\begin{equation}\label{LiapunovA}
\sum_{k=-\infty}^m (d_{n,k,t}^{(j)})^2\left (\sum_{k=-\infty}^n (d_{n,k,1}^{(j)})^2\right )^{-1}\le C(t).
\end{equation}
Then it is easy to see that (\ref{Liapunov}) will follow if
\begin{equation}\label{Liapunov1}
{\tilde L}^{(j)}(\d,n, t)=\left ((A_n^{(j)})^{-1}\max_{-\infty<k\le m} |d_{n,k,t}^{(j)}|\right )^{\d} \to 0.
\end{equation}
This quantity with $\d=1$ was estimated in the case of non-tapered filter in \cite{Paul21}, the estimation in the case of tapered filter goes along the same lines as in \cite{Paul21}, therefore we shall provide the proof of (\ref{Liapunov1}) only for several $j$. Let us denote
\begin{equation}\label{Liapunov1A}
I_1^{(j)}=\max_{0\le k <\infty}\big |\sum_{i=1}^{m}{\bar a}_{i+k}\big |, \quad I_2^{(j)}=\max_{0<k\le m}\big |\sum_{i=k}^{m}{\bar a}_{i-k}\big |,
\end{equation}
then (\ref{Liapunov1}) will follow if we show
\begin{equation}\label{Liapunov2}
\left (A_n^{(j)} \right )^{-1}\max (I_1^{(j)}, I_2^{(j)}) \to 0.
\end{equation}
Let us consider the case $0<\g_1<1$, then, for any fixed $t$, $m_1=[n^{\g_1}]< m$ and we have, for $j=1, 2, 3$,
$$
I_1^{(j)}=\max_{0\le k <\infty}\Big |\sum_{i=1}^{m}a_{i+k} \ind{\{0\le i+k \leq m_1\}}\Big |=\max_{0\le k \le m_1-1}\Big |\sum_{i=1}^{ m_1-k}a_{i+k}\Big |\le \sum_{i=1}^{m_1}|a_{i}|
$$
In a similar way we get
$$
I_2^{(j)}=\max_{0\le k \le m}\big |\sum_{i=k}^{m}{\bar a}_{i-k}\big |=\max_{0\le k \le m}|\sum_{i=0}^{\min(m, m_1-k)}a_{i}|\le \sum_{i=0}^{ m_1}|a_{i}|.
$$
From these two estimates we get
$\max (I_1^{(1)}, I_2^{(1)})\le C m_1^{1-\b} , \ \max (I_1^{(j)}, I_2^{(j)})\le C, \ {\rm for} \ j=2, 3.$
Taking into account the expressions of $A_n^{(j)} $ from Proposition \ref{prop1} we get (\ref{Liapunov2})for $j=1, 2, 3.$

In the case $\g_1>1$ we have $m_1>m$ and in a similar way we can get the following estimates
$$
I_1^{(j)}=\max_{0\le k \le m_1-1}\Big |\sum_{i=1}^{\min (m, m_1-k)}a_{i+k}\Big |\le \max \left (\sum_{i=1}^{ m}|a_{i}|, \sum_{i=m_1-m}^{ m_1}|a_{i}| \right )
$$
and
$$
I_2^{(j)}\le \max_{0\le k \le m}\Big |\sum_{i=0}^{\min (m_1, m-k)}a_{i}\Big |\le \sum_{i=0}^{ m}|a_{i}|.
$$
Again, from these two estimations and the expressions of $A_n^{(j)}$ we get (\ref{Liapunov2})for $j=4, 5, 6.$

In the case $\g_1=1$ and $m_1=[nc]$ we must consider two cases $m<m_1$ and $m\ge  m_1$ and in the same way as above we  get (\ref{Liapunov2}) for $j=7, 8, 9.$ Thus, we had proved (\ref{Liapunov}).
\halmos
Propositions \ref{prop1} and \ref{prop2} prove  \ref{fddconv}.

As it was mentioned in Remark \ref{rem2}, main tool in the proof of (\ref{wcconv}) is Proposition \ref{prop1}. The tightness can be proved in the same way as it is done in \cite{Sabzikar}, using the ideas from \cite{Billingsley}. Namely, we use the tightness criterion given in  (6.9) in \cite{Sabzikar} with $p=2+\d$ and we need to estimate the quantity $E|Z_n(t)-Z_n(s|^p)$ for $0<s<t\le 1$. Let us introduce notations $m_n(t)=[nt]/n, \ m(n, t,s)=[nt]-[ns]$. Using the stationarity of a linear process and equality $m(n, t,s)=[n(m_n(t)-m_n(s))]$ we have
$$
E|Z_n(t)-Z_n(s|^p)=E|Z_n (m_n(t)-m_n(s))|^p.
$$
By Rosenthal's inequality (see, for example \cite{Giraitis}, Prop. 4.4.3) we arrive at estimation of $\left (E|Z_n (m_n(t)-m_n(s))|^2\right )^{p/2}$, and this can be done by means of Proposition \ref{prop1}. For all cases $j=1, 2,\dots, 9$ we must verify that $H(j)p>1$ (for $j=7, 9$ we must use asymptotic of $C_{10}(t, \b, c)$ given in (\ref{VarZn5}), similar asymptotic is for $C_{17}(t, \b, c)$). If $H(j)=1/2$, then this inequality holds since $p>2$. If $H(j)=3/2-\b$ and $1/2<\b<1$, then again $H(j)p>1$ if $\d>0$ and only in the case of ND, where $1<\b <3/2$, we need stronger moment condition $p>2 (3-2\b)^{-1}$. (\ref{wcconv}) and at the same Theorem \ref{thm1} is proved.
\halmos

\subsection{The stable case}\label{subsec2}

In this subsection we consider limit theorems for linear processes with tapered filter (\ref{newtaper}) but with heavy tailed innovations. Namely, let
\begin{equation}\label{linpr2}
X^{(n)}=\{X_k^{(n)} \}, \quad   X_k^{(n)}=\sum_{i=0}^\infty {\bar a}_i\eta_{k-i}, \ k\in \bz,
\end{equation}
where ${\bar a}_i$ are defined in (\ref{newtaper}) and $\{\eta_i,\ i\in \bz\}$ are heavy tailed i.i.d. random variables. Instead of taking the most general case of innovations, belonging to the Domain of attraction of  an $\a$-stable law, as it is in the pioneering paper \cite{Astrauskas}, we shall consider particular case of innovations with Pareto distribution. This was done in our previous papers \cite{Paul21}, \cite{Paul22}, and \cite{PaulDam5}, and we can explain the reason for such choice.  Domain of attraction of a Gaussian law is very broad - all random variables with finite variance belong to this domain. Contrary, domain of attraction of a particular  $\a$-stable law is very narrow, only random variables having  tails  of the form $x^{-\a}L(x)$ with  this particular $\a$, belong to it. (for strict definitions see, for example , \cite{Ibragimov}). Here $L$ is a slowly varying function, which  has no effect to a limit stable distribution, only the normalizing sequences depend on this function, see \cite{Astrauskas}. On the other hand, in the proofs of limit theorems in \cite{Astrauskas} main difficulties arise dealing with slowly varying functions, especially that another slowly varying function is assumed describing the asymptotic of the coefficients $a_i$. In \cite{Damarackas} even it was noted that one statement in \cite{Astrauskas}, connected with slowly varying functions, is incorrect (see Remark 1 in \cite{Damarackas}; fortunately, the proof can be corrected as shown in \cite{Damarackas}). Another fact, which motivated us to choose Pareto distribution instead of a general distribution in the domain of attraction of a stable law, was described in our old paper \cite{PaulJuoz}. In this paper it was shown, that if we add logarithmic function to Paretian tail, e.g., consider density of the form $p(x)=Cx^{-\a-1}|\ln x|^\g, \ x>x_0$, then we must take normalizing constant (for a centered sum of $n$ i.i.d. random variables with this density) of  the form $Cn^{1/\a}L_1(x)$ with another slowly varying function $L_1$ and we get only logarithmic rate of convergence (with respect to $n$) to a limit stable law. Moreover, even an asymptotic expansion in this example is possible only in negative  powers of $\ln n$. This shows the big difference between domain of normal attraction (when the normalizing constants are $Cn^{1/\a}$ and  $L(x)=C+o(1)$) and  domain of attraction of a stable law (when the normalizing constants are $Cn^{1/\a}L_1(x)$ with $L_1(x)\to \infty\ {\rm or} \ 0$ as $x\to \infty$). 
If one looks at log-log plots of real data from various fields in \cite{Kagan}, \cite{Aban}, and \cite{Meerschaert} one can see the tapering effect, therefore during last decades there were papers dealing with estimation of parameters of  tapered (including truncated) Pareto distributions, see the above cited papers \cite{Aban} and \cite{Meerschaert}. Tapered Pareto distribution, being the most simple between tapered other heavy-tailed distributions, quite well represents these other distributions, for example, in \cite{Meerschaert} there was generated  a sample from tapered  stable distribution and it turned out that tapered Pareto model gave quite good fit to this data.
Taking into account that in practice even estimation of the tail index $\a$ is quite difficult problem, not speaking about statistical procedures which would allow to say what a slowly varying function is present in the tail behavior of a distribution under consideration, it is natural to concentrate our studies of limit theorems for Paretian distributions.

Let $\theta =\theta (\a)$ stand for  the standard Pareto distribution with probability density and distribution functions
 \begin{equation}\label{paretodf}
f(x)= \a  x^{-\a-1}, \quad F(x)=1-x^{-\a}, \ x\ge 1,
\end{equation}
where $\a>0$. Let us denote $\eta=\theta$, if $0<\a<1$ and $\eta=\theta -E\theta$, if $1<\a<2$. Let $\eta_i, \ i\in \bz$ be i.i.d. copies of $\eta$. We exclude the case $\a=1$, since usually in this case, in order to avoid problems with centering, it is assumed symmetric distribution, while we on purpose had chosen the case of completely asymmetrical case of positive Pareto distribution.
In this subsection we consider $S_{n}(t; X^{(n)})$ where $S_{n}(t; X)$ is defined in (\ref{sum}) and $ X^{(n)}$ - in (\ref{linpr2}). As in subsection 2.1 we use the notation $Z_n(t)=A_n^{-1} S_n(t, X^{(n)} ),$ where $A_n$ is a normalizing sequence for $ S_n(1, X^{(n)})$.
 For a filter $\{a_i\}$ we introduce three sets of conditions, in the cases of LRD and ND  assuming condition (\ref{cond1}):

(i) $1/\a<\b<1$ - the case of LRD;

 (ii) $\sum_{i=0}^\infty |a_i|<\infty$ and  $\sum_{i=0}^\infty a_i\ne 0$ - the case of SRD;

 (iii) $\max (1, 1/\a)<\b<1+1/\a$ and $\sum_{i=1}^\infty a_i=0$, the case of ND.

As in the Gaussian case, combining three cases of tapering with three cases of dependence,  we  consider nine cases, and we shall number these nine cases in the same order as in subsection 2.1.
Then we shall add index $j$ to most quantities, introduced above, but we omit this index, if it will be clear what case is considered.

Let $ M_\a$ be an $\a$-stable random measure on $\bbr$ with the Lebesgue control measure and constant skewness intensity ${\bar \b}=-1$. For $j=1, 2, \dots, 9$ we define the following processes
\begin{equation}\label{Udef}
U_{j}(t)=\int_{\RR}\H^{(j)}_{\b}(u,t,c) M_\a(\du),
\end{equation}
where
\begin{equation}\label{Hdef}
\H_{\b}^{(j)}(u,t,c)=\begin{cases}
\sum_{k=0}^{\infty}a_k\ind{[0,t)}(u), &\text{ if } j=1,2,3,5,8,\\
\left ((t-u)_+^{1-\b}-(-u)_+^{1-\b}\right ) (1-\b)^{-1}, &\text{ if } j=4,6, \\
\frac{\ind{[-c, t]}(u)}{1-\b}\left (\left (\min (t-u, c)\right )^{1-\b}-(-u)_+^{1-\b}\right ), &\text{ if } j=7,9.\\
\end{cases}
\end{equation}

\begin{teo}\label{thm3} For $Z_n(t)$, defined by means of a linear random process (\ref{linpr2}) with the tapered  filter, for all  $j=1, 2, \dots, 9$,
the following relation holds
\begin{equation}\label{main}
\left \{Z_{n}^{(j)}(t), \ t \in \bbr_+\right \}\fdd \{U_{j}(t),  \ t \in \bbr_+ \}.
\end{equation}
\end{teo}

We see that in all cases, except cases $j=7, 9$ we get the same limit processes which  are in the cases of non-tapered filter in \cite{Astrauskas}, only in the cases $j=7,9$ (the case of moderate tapering and LRD or ND) we get new limit processes.

{\it Proof of Theorem \ref{thm3}}. In the proof we use the scheme of the proof of Theorem 2 in \cite{PaulDam5}, where the important result of Astrauskas \cite{Astrauskas} was generalized to random linear fields with factorizing coefficients of the filter. Although we deal with linear processes, for us it is more convenient to use \cite{PaulDam5}, since in \cite{Astrauskas}, not all cases of dependence are proved and only in symmetric case of innovations, while we deal with completely asymmetric case. Moreover, now we consider particular case of innovations ($h(x)\equiv 1$ in (2.13) in \cite{PaulDam5}), therefore the proof is more simple comparing with the proof in \cite{PaulDam5}. Since the beginning of the proof is the same for all cases  $j=1, 2, \dots, 9$, at the beginning we shall skip this index. In order to prove the convergence of finite-dimensional distributions we must consider linear combination $\sum_{l=1}^m x_lZ_n(t_l)$ with an arbitrary fixed $m$. It is clear that it sufficient to consider the case $m=1, x_1=1, t_1= t.$ Let us consider the quantity
$$
 J_n :=E \exp\left( \ri  Z_{n}(t) \right)= E \exp\left( \ri  A_{n}^{-1} \sum_{i=-\infty}^{\sv{nt}}d_{n, i, t}\eta_{i} \right)
$$
where  $d_{n,i,t}=\sum_{k=1}^{[nt]}{\bar a}_{k-i}$ for $i\le 0$ and $d_{n,i,t}=\sum_{k=i}^{[nt]}{\bar a}_{k-i}$ for $i> 0$. Therefore, using formulae (2.13) and (2.14) from \cite{PaulDam5}, only with $c_1=0$ and $h(x)\equiv 1$, we get
\begin{eqnarray}\label{lnJn}
\ln  J_n &=& -\sum_{i=-\infty}^{\sv{nt}}\left( 1+r\left( A_{n}^{-1}d_{n, i, t}  \right) \right)\abs{ A_{n}^{-1}d_{n, i, t} }^\a \\ \nonumber
    &+& \ri{\bar \b}\tau_\a\sum_{i=-\infty}^{\sv{nt}}\left( 1+r\left( A_{n}^{-1}d_{n, i, t}  \right)  \right)(A_{n}^{-1}d_{n, i, t})^{<\a>}.
\end{eqnarray}
Here  $r(t)\rightarrow 0,$ as $ t\rightarrow 0$, ${\bar \b}=-1$, $\tau_\a=\tan\left( {\pi\a}/{2} \right)$ if $\a\neq 1$, and $x^{\langle\a\rangle}=\abs{x}^\a\sign (x)$. Since the functions $|x|^\a$ and $x^{<\a>}$ in many respects are similar, it is clear that, denoting by $f(x)$ one of these functions, we must investigate the asymptotic behavior of the quantity
\begin{equation}\label{In}
 I_n=\sum_{i=-\infty}^{\sv{nt}}f\left( A_{n}^{-1}d_{n, i, t}  \right) \left( 1+r\left( A_{n}^{-1}d_{n, i, t}  \right) \right).
\end{equation}
In all cases we shall prove that $ A_{n}^{-1}d_{n, i, t} \to 0$ as $n\to \infty$, uniformly with respect to $i\in \bz$, therefore we must investigate the remaining sum in (\ref{In}). As in \cite{PaulDam5}, we write this sum as integral
\begin{eqnarray}\label{In1}
{\tilde I}_n &:=& \sum_{i=-\infty}^{\sv{nt}}f\left( A_{n}^{-1}d_{n, i, t}  \right)=\int_{-\infty}^{\sv{nt}}f\left( A_{n}^{-1}d_{n, \sv{u}, t}  \right)\dd u  \\ \nonumber
   &=& \int_{-\infty}^{t}n f\left( A_{n}^{-1}d_{n, \sv{nu}, t}  \right)\dd u.
\end{eqnarray}
Now we must prove the point-wise convergence (for a fixed $t$)
\begin{equation}\label{bound}
\frac{d_{n, \sv{nu}, t} }{z_{\b,\g_1, n}}\rightarrow \H_{\b}(u,t)
\end{equation}
with some sequence $z_{\b, \g_1, n}$ and some function $\H_\b$. Also, in order to apply the dominated convergence theorem, we need the bound
\begin{equation}\label{bound1}
\frac{|d_{n, \sv{nu}, t} |}{z_{\b, \g_1, n}}\leq  G_{\b}(u,t),
\end{equation}
where $G_{\b}(u,t)$ is a function, for a fixed $t$ satisfying
\begin{equation*}\label{integrableBound}
\int_{-\infty}^{\infty} \abs{G_{\b}(u,t)}^\a \dd u<\infty.
\end{equation*}
Relations (\ref{bound}) and (\ref{bound1}) for quantities $d_{n, i, t}$, expressed via initial filter coefficients $\{ a_i\}$ in all cases of dependence were proved in \cite{PaulDam3}, Prop. 4, and used in \cite{PaulDam5}, Prop. 1. Now we must prove the same relations (\ref{bound}) and (\ref{bound1}) in all nine cases $j=1, 2, \dots , 9$ for  $d_{n, i, t}$, expressed via tapered filter
$\{{\bar a}_i\}$.
We start with the case $j=1$,
 e.g., $0<\g_1<1$ and $1/\a<\b<1$. Writing
$$
d_{n, \sv{nu}, t}^{(1)} =\sum_{k=\left( -\sv{nu} \right)\vee 0 }^{\sv{nt}-\sv{nu}}{\bar a}_{k}=\sum_{k=\left( -\sv{nu} \right)\vee 0 }^{\sv{nt}-\sv{nu}} a_{k}\ind{[ 0\leq k \leq n^{\g_1}]},
$$
it is easy to see that for $u<0$ we have
$d_{n, \sv{nu}, t}^{(1)} \to 0, \quad {\rm as} \ n\to \infty,$
since for any fixed $u<0$ for sufficiently large $n$ we have $-\sv{nu}>n^{\g_1}$ and the above written sum is empty.
For $0<u<t$ we have
$$
d_{n, \sv{nu}, t}^{(1)}=\sum_{k= 0 }^{\sv{nt}-\sv{nu}} a_{k}\ind{[ 0\leq k \leq n^{\g_1}]}=a_0+\sum_{k=1}^{n^{\g_1}}k^{-\b}.
$$
therefore, taking $z_{\b, \g_1, n}^{(1)}=n^{\g_1(1-\b)}$, we get
$d_{n, \sv{nu}, t}^{(1)}(z_{\b, \g_1, n}^{(1)})^{-1}\rightarrow  (1-\b)^{-1}\ind{[0,t)}(u).$
To get (\ref{bound1}), let us note that $-\sv{nu}>n>n^{\g_1}$ for $u<-1$, therefore we take $G_{\b}^{(1)}(u,t)\equiv 0$ for $u\in (-\infty, -1)\cup (t+1, \infty)$. It is easy to see that for $u\in [-1, t+1]$ we can take $G_{\b}^{(1)}(u,t)=|a_0|+(1-\b)^{-1}$.

The case $j=3, 0<\g_1<1$, $\max (1, 1/\a)<\b<1+1/\a$, and $\sum_{i=0}^\infty a_i=0$ is similar to the case $j=1$, and for this reason  the proof was not provided nor in \cite{PaulDam4} neither in \cite{PaulDam5}. Therefore, we shall give the sketch of the proof. Now we use equality $\sum_{i=0}^n a_i=-\sum_{i=n+1}^\infty a_i$ and we will adopt the following definition for integer part of negative numbers: for $a<0$,  $\sv{a}=\min (n\in \bz: n\ge a)$.  For $-1/n<u<t$ we easily get
$$
d_{n, \sv{nu}, t}^{(3)} =\sum_{k=\left( -\sv{nu} \right)\vee 0 }^{\sv{nt}-\sv{nu}}{\bar a}_{k}=\sum_{k=0 }^{n^{\g_1}} a_{k}=-\sum_{k=\sv{n^{\g_1}}+1}^\infty a_{k},
$$
therefore, taking $z_{\b, \g_1, n}^{(3)}=n^{\g_1(1-\b)}$, we  get
$$
d_{n, \sv{nu}, t}^{(3)}(z_{\b, \g_1, n}^{(3)})^{-1}\rightarrow  -(1-\b)^{-1}\ind{[0,t)}(u).
$$
It is easy to see that if $u<-(n^{\g_1}+1)/n$, then $-\sv{nu}>\sv{n^{\g_1}}$ and $d_{n, -\sv{nu}, t}^{(3)}=0$. Since $\g_1<1$ then for any fixed $u<0$ for sufficiently large $n$ we shall have $u<-(n^{\g_1}+1)/n$, therefore, for $u<0$, we have
$d_{n, -\sv{nu}, t}^{(3)}(z_{\b, \g_1, n}^{(3)})^{-1}\rightarrow  0.$
Since $d_{n, \sv{nu}, t}^{(3)} =0$ for all $n\ge 1$ if $u\in (-\infty, -1)\cup ((t+1, \infty)$, it is  easy to construct the majoring function
$$
G_{\b}^{(3)}(u,t)=\ind{[-1, 0)}(u)\left ( (t-u+1)^{1-\b} -(-u)^{1-\b}\right )(1-\b)^{-1}+\ind{[0, t]}(u)(1-\b)^{-1}.
$$
The case $j=2$ is easy, taking $z_{\b, \g_1, n}^{(2)}\equiv 1$ we get the relation and the majoring function
$$
d_{n, \sv{nu}, t}^{(2)}\rightarrow  \ind{[0,t)}(u)\sum_{i=0}^\infty a_i, \quad G_{\b}^{(2)}(u,t)=\ind{[-1, t]}\sum_{i=0}^\infty |a_i|.
$$

Now let us consider cases $j=4, 5, 6$ with weak tapering $\g_1>1$. It is easy to see that for any fixed $u$ and $t$, for sufficiently large $n$,  $\sv{nt}-\sv{nu}<n^{\g_1}$ and $ -\sv{nu}<n^{\g_1}$, therefore for any fixed $u, t$,  sufficiently large $n$, and for all $j=4, 5, 6$ we have
$$
d_{n, \sv{nu}, t}^{(j)}=\sum_{k=\left( -\sv{nu} \right)\vee 0 }^{\sv{nt}-\sv{nu}}{\bar a}_{k}=\sum_{k=\left( -\sv{nu} \right)\vee 0 }^{\sv{nt}-\sv{nu}}a_{k}.
$$
Thus, we have the situation with non-tapered filter, which was considered in \cite{PaulDam4}, therefore we obtain (\ref{bound}) for all $j=4, 5, 6$ and (\ref{bound1}) for $j=4, 6$. There is a problem with a bound in the case $j=5$, since in \cite{PaulDam4} this case was proved under a little bit stronger condition - instead of condition $\sum_{i=0}^\infty |a_i|<\infty$ it was required $\b>1$. At present we do not know how to estimate $|d_{n, \sv{nu}, t}^{(5)}|$ for $u<0$, having only condition $\sum_{i=0}^\infty |a_i|<\infty$. On the other hand, in \cite{Astrauskas} the limit theorem is proved for non-tapered filter under the condition $\sum_{i=0}^\infty |a_i|<\infty$, but the scheme of the proof is quite different. The idea is to take $m=m(n)\to \infty$ and $m=o(n)$ and to write
$$
Z_n (t)=A_n^{-1}\sum_{i=0}^m a_i \sum_{k=1}^{\sv{nt}}\eta_k +W_n(t).
$$
The first term in this decomposition gives us the limit process, and it is shown that $W_n(t)\pc 0$, as $n\to \infty$. But in the proof of this step  in \cite{Astrauskas}  truncation of innovations $\eta_k$ with different truncation level for each $k$ is used. We had verified this proof in our setting, and it works, but since it takes additional 3-4 pages of text we decided to leave in formulation of Theorem \ref{thm3}  assumption $\sum_{i=0}^\infty |a_i|<\infty$, while proof which we provide in the paper for (\ref{bound1}) is under condition $\b>1$.

It remains the last three cases $j=7, 8, 9$ of moderate tapering, $\g_1=1, \l(n)=cn$. We start with the case $j=7$, and consider
$$
d_{n, \sv{nu}, t}^{(7)} =\sum_{k=\left( -\sv{nu} \right)\vee 0 }^{\sv{nt}-\sv{nu}} a_{k}\ind{[0\le k \le cn]}.
$$
The location of $t$ with respect to $c$ is important, at first let us consider the case $c<t$. For $0\le u<t$, taking $z_{\b, n}^{(7)}=n^{1-\b}$ it is easy to get the following relation
$$
d_{n, \sv{nu}, t}^{(7)}(z_{\b, n}^{(7)})^{-1}\rightarrow (1-\b)^{-1}\left (\min (t-u, c)\right )^{1-\b}, \quad  {\rm as} \ n\to \infty.
$$
For $u<0$ let us note that if $u<-c$, then for sufficiently large $n$ we have $-\sv{nu}>\sv{cn}$ (for $u<-c-1$ this inequality holds for all $n$) and therefore $d_{n, \sv{nu}, t}^{(7)}=0$. For $-c<u<0$ we have
$$
d_{n, \sv{nu}, t}^{(7)}=\sum_{k= -\sv{nu} }^{\sv{nt}-\sv{nu}} a_{k}\ind{[0\le k \le cn]}\ind {[-\sv{nu}<\sv{cn}]},
$$
therefore, taking into account that for $u<0$ and $c<t$ we have $\min (t-u, c)=c$,   it is not difficult to obtain the following relation
$$
d_{n, \sv{nu}, t}^{(7)}(z_{\b, n}^{(7)})^{-1}\rightarrow (1-\b)^{-1}\left ( c^{1-\b}- (-u)^{1-\b}\right ), \quad  {\rm as} \ n\to \infty.
$$
In the case $c\ge t$ in a similar way we get the following results. For $0\le u<t$ we have $\min (t-u, c)=t-u$ and
$$
d_{n, \sv{nu}, t}^{(7)}(z_{\b, n}^{(7)})^{-1}\rightarrow (1-\b)^{-1}(t-u)^{1-\b}, \quad  {\rm as} \ n\to \infty.
$$
If $u<0$, then again we need to consider only interval $-c<u<0$ and we get
$$
d_{n, \sv{nu}, t}^{(7)}(z_{\b, n}^{(7)})^{-1}\rightarrow (1-\b)^{-1}\left (\min (t-u, c)\right )^{1-\b}\ind{[-c, 0)}(u), \quad  {\rm as} \ n\to \infty.
$$
It is easy to see that all these cases can be written in one expression as given in (\ref{Hdef}). Since  $d_{n, \sv{nu}, t}^{(7)}=0$ outside of the interval $-c-1<u<t+1$, there is no difficulties to construct the function $G_{\b}^{(7)}(u,t)$, since $1-\b>0$.

The case $j=8$ is simple. Taking $z_{\b, n}^{(8)}\equiv 1$ it is easy to get
$$
\frac{d_{n, \sv{nu}, t}^{(8)}}{z_{\b, n}^{(7)}}\rightarrow \sum_{k=0}^{\infty}a_k\ind{[0,t)}(u) \quad  {\rm as} \ n\to \infty.
$$
Since outside of the interval $-c-1<u<t+1$ function $d_{n, \sv{nu}, t}^{(8)}\equiv 0$, thus we can take $G_{\b}^{(8)}(u,t)=\sum_{k=0}^{\infty}|a_k|\ind{[-c-1, t+1]}(u)$.

The case $j=9$ is similar to the case $j=7$ (see also the sketch of the proof in the case $j=3$). Taking $z_{\b, n}^{(9)}=n^{1-\b}$   we get the same expression of limit function $\H^{(9)}$ as in the case $j=7$. But here it is necessary to note that functions $\H^{(7)}$ and $\H^{(9)}$ differs essentially - if $\H^{(7)}$ is a bounded function, due to the condition $\b<1$, function $\H^{(9)}$ is unbounded (tends to infinity when $u\to 0$ or $u\to t$), but $|\H^{(9)}|^\a$  is integrable  on  interval $(-c, t)$ due to the condition $\b<1+1/\a$. Again,  outside of the interval $-c-1<u<t+1$ function $d_{n, \sv{nu}, t}^{(9)}\equiv 0$, therefore it is not difficult to construct the function $G_{\b}^{(9)}(u,t,c)$. One can take
$$
G_{\b}^{(9)}(u,t,c)=\frac{(t-u)^{(1-\b)}+(-u)^{1-\b}}{\b-1}\ind{(-c-1, 0)}(u)+ \frac{(\min (t-u, c))^{1-\b}}{\b-1}\ind{[0, t)}(u).
$$
Thus, for all cases  $j=1, 2, \dots, 9$ we have proved  relations (\ref{bound}) and (\ref{bound1}). Having these relations, for all $j=1, 2, \dots, 9$, we easily get the relation
\begin{equation}\label{bound2}
(A_{n}^{(j)})^{-1}d_{n, i, t}^{(j)} \to 0,
\end{equation}
as $n\to \infty$, uniformly with respect to $i\in \bz$. Namely, writing (with $i=\sv{nu}$)
$$
\frac{|d_{n, i, t}^{(j)}|}{A_{n}^{(j)}}=\frac{1}{n^{1/\a}}\frac{|d_{n, \sv{nu}, t}^{(j)}|}{z_{\b, n}^{(j)}}
$$
and in all cases $j\ne 6, 9$, where functions $G_{\b}^{(j)}$ from (\ref{bound1}) are bounded, we get (\ref{bound2}). Functions $G_{\b}^{(j)}, \ j=6, 9$ are unbounded (there is a term $(-u)^{1-\b}\ind{(-c-1, 0)}(u)$) with $\b>1$, but then we can simply estimate $|d_{n, \sv{nu}, t}^{(j)}|\le \sum_{k=0}^\infty |a_k|$ and we  get (\ref{bound2}) (we recall that $A_{n}^{(j)}=n^{1-\b+1/\a}$ and $\b<1+1/\a$).

Now we can finish the proof of the theorem. Collecting (\ref{lnJn})-(\ref{bound}),  (\ref{bound2}) and applying the dominated convergence theorem we get the relation (\ref{main}).
\halmos

\subsection{New limit processes and some non-standard situations}

As it was mentioned above, we get different limit processes (comparing with limit processes in the case of non-tapered filters) only in the case of moderate tapering and cases of LRD and ND. The same result was in \cite{Sabzikar}, but limit processes in this paper and for these two cases differ from our limit processes. This difference is caused by different tapering of the filter: exponential tapering  in  \cite{Sabzikar} and tapering by indicator function (or truncation) in our paper.
In \cite{Sabzikar} and \cite{Sabzikar2} the tempered fractional stable motion of the second kind for $1<\a\le 2$ and abbreviated as TFSMII was introduced
\begin{equation}\label{TFSMII}
Z_{H, \a, \l}^{II}(t)=\int_{-\infty}^{\infty} h_{H, \l} (t; u) M_{\a}(du),
\end{equation}
where $H=d+1/\a>0, \ \l> 0$, $M_\a$ is an independently scattered stable random measure with the exponent $\a$ and the Lebesgue control  measure, and
\begin{eqnarray*}
h_{H, \l} (t; u) &=& (t-u)_+^{H-\frac{1}{\a}}\exp (-\l (t-u)_+)- ((-u)_+)^{H-\frac{1}{\a}}\exp (-\l (-u)_+) \\
   & + & \int_0^t (s-u)_+^{H-\frac{1}{\a}}\exp (-\l (s-u)_+) ds.
\end{eqnarray*}
 In the case $\a=2$ we have tempered fractional Brownian motion of the second kind (TFBMII), in the names of these processes there are the words "second kind", since in \cite{Sabzikar1} and \cite{Sabzikar3} there were introduced similar, but different, TFBM and TFSM, respectively.
%
Processes $U_j, \ j=7,9$ were defined in (\ref{Udef}) by means of an $\a$-stable random measure on $\bbr$ with the Lebesgue control measure and constant skewness intensity ${\bar \b}=-1$, since we consider one-sided Pareto random variables. Clearly, in (\ref{Udef}) we can take $M_\a$ with any skewness intensity ${\bar \b}\in [-1, 1]$, so we propose the following definition.
 For $0<\a<2$ let us denote by $M_{\a,{\bar \b}},\s$ an $\a$-stable random measure on $\bbr$ with   control measure $\mu(dx)=\s^{1/\a}\l(dx)$, where $\l$ is Lebesgue measure,   and constant skewness intensity ${\bar \b}\in [-1, 1]$. For $\a=1$ we assume ${\bar \b}=0$. If $\a=2$, then  $M_{2,{\bar \b}, \s}:=M_{2, \s}$ stands for random Gaussian measure with   control measure $\s^{1/2}\l(dx)$. If $\s=1$, we skip this parameter from notations.
 \begin{definition}\label{TFSMIII} Let $0<\a\le 2, H>0, c>0$. The random process
 \begin{equation}\label{defTFSMIII}
Z_{H, \a, c, \s}^{III}(t)=\int_{-\infty}^{\infty} {\tilde h}_{H, c} (t; u) M_{\a, {\bar \b}, \s}(du),
\end{equation}
where
\begin{equation}\label{deftildeh}
{\tilde h}_{H, \a, c} (t; u)=(t-u)^{H-\frac{1}{\a}}\ind{(-c, t)}(u)- ((-u))^{H-\frac{1}{\a}}\ind{(-c, 0)}(u),
\end{equation}
is called tapered fractional stable motion of the third kind and abbreviated as TFSMIII. In case $\a=2$ this process is called tapered fractional Brownian motion of the third kind (TFBMIII) and is denoted by $B_{H, c, \s}^{III}(t)$.
\end{definition}
Since the function ${\tilde h}_{H, \a, c} (t; u)$ is zero outside of interval $(-c, t)$, the existence of stochastic integral in (\ref{defTFSMIII}) follows easily  from the fact, that only point $u=0$ must be verified: $(H-1/\a)\a >-1$ for all $H>0$. This is effect of tapering function (in our case it is indicator function), for FSM we need the condition $0<H<1$. Definition of TFSMIII in (\ref{defTFSMIII}) is similar to definition of TFSM in \cite{Sabzikar3}, changing the exponential taper by indicator function. By the way, both words "tempered" and "tapered", which can be considered as synonyms, give the same letter T in abbreviation.

To compare limit processes $U_j, j=7,9$ (which now can be called TFSMIII) with limit processes from (\ref{TFSMII}) in the stable case it is easy. Correspondence between parameters, used in \cite{Sabzikar}, and ours is  $\ d=1-\b,$ for $d\ne 0 $ and $\l_*=c^{-1}$, therefore the exponent $H-1/\a=1-\b$  in the expression of $h_{H, \l_*} (t; u)$ is the same as in $\H_{\b}^{(j)}(u,t,c), \ j=7,9$. But there is  main difference between these two functions: the support of the functions $\H_{\b}^{(j)}(u,t,c), \ j=7,9$ is a finite interval $[-c, t]$, while function $h_{H, \l_*} (t; u)$ is supported on $(-\infty, t]$ with exponential decay as $u\to -\infty$, moreover, (\ref{TFSMII}) is defined only for $1<\a\le 2$.
 More difficult case is Gaussian, since  all three types of TFMB, introduced in \cite{Sabzikar1}, (\ref{TFSMII}), and (\ref{defTFSMIII}), are  expressed as stochastic integrals, while processes    $U^{j}, j=7, 9$ in Theorem \ref{thm1} are defined by means of variances and covariances. For comparison we proceed as follows. We investigate the asymptotic behavior as $t\to 0$ and $t\to \infty$ of  variance $W^{(7)}(t)$, given in (\ref{VarZn3}), then we investigate asymptotic of variance of $B_{H, c}^{III}(t)$.
Looking at the expressions of variances $W^{(j)}(t), \ j=7, 9$ it is easy to see that expressions of functions $C_{10}$ and $C_{17}$ are functions on $z=c/t$, integrands in these functions are power functions, so all integrals, except one, can be integrated and we get sums of powers of $z$ and  $z-1$. Only integrating integral in the expression of $C_{1}$ we shall  get the integral $I(z)=\int_0^{z-1} (y(1+y))^{1-\b} dy$, for $z>1$. But in \cite{Paul21}, see (2.6) and (2.7) therein, it was proved the finiteness of the integral $C_0(\b)=:\int_0^{\infty}\left (\int_0^{1}(x+y)^{-\b} dx\right )^2 dy$ for $1/2<\b<3/2, \b\ne 1$, therefore evaluating $C_1$ we write
$$
\int_0^{z-1}\left (\int_0^{1}(x+y)^{-\b} dx\right )^2 dy=C_0(\b)- \int_{z-1}^{\infty}\left (\int_0^{1}(x+y)^{-\b} dx\right )^2 dy
$$
and for integrand and large $y$ we use the relation $\left ((1+y)^{1-\b}-y^{1-\b} \right )^2= (1-\b)y^{-2\b}+O(y^{-2\b-1})$ and the fact that $2\b>1$.
Let us consider the asymptotic of $W^{(7)}(t)$, given in (\ref{VarZn3}), for $t\to 0$ and for $t\to \infty$.
To get asymptotic behavior for $t \to 0$ we consider $C_{4}(t, \b, c)$, while $C_{8}(t, \b, c)$ gives asymptotic  behavior for $t \to \infty$  ($C_{9}(1, \b, c)$ does not depend on $t$, it affects only constants in the asymptotic).
 Performing some calculations we can get that for $0<t<c$ and as $t\to 0$,
\begin{equation}\label{VarZn5}
\lim_{n\to \infty}{\rm Var}Z_n^{(7)}(t)=C_9^{-1}(1, \b, c)\left (C_0(\b)+C_3(\b)\right )t^{3-2\b}+O(t^2).
\end{equation}
Here $C_0(\b)=:\int_0^{\infty}\left (\int_0^{1}(x+y)^{-\b} dx\right )^2 dy$, and finiteness of this integral for $1/2<\b<3/2, \b\ne 1$ was shown in \cite{Paul21}, see (2.6) and (2.7) therein.
In a similar way we can investigate $C_{8}(t, \b, c)$, as $t\to \infty$, and to find that   as $t \to \infty$,
\begin{equation}\label{VarZn6}
\lim_{n\to \infty}{\rm Var}Z_n^{(7)}(t)=C_9^{-1}(1, \b, c)\left (\frac{c^{2(1-\b)}}{(1-\b)^{2}}t -\frac{2}{(1-\b)(2-\b)(3-2\b)}c^{3-2\b}\right ).
\end{equation}
 We see  that $W^{(7)}(t)$ for $t\to 0$  has the same asymptotic behavior  as variance of FBM with Hurst parameter $H(7)$  and for $t\to \infty$ - as variance of BM.

Since the integrand function, defined in (\ref{deftildeh}) for $\a=2$, is expressed as power functions, it is not difficult to calculate variance of  $B_{H, c}^{III}(t)$. In order to compare this variance with $W^{(7)}(t)$ we take  $\s=1$ and $H-\frac{1}{2}=1-\b, \ 1/2 <\b<1$  in (\ref{deftildeh}) and we evaluate ${\tilde W}_{\b,c}(t):=E\left (B_{1-\b, c}^{III}(t)\right )^2$,
\begin{equation*}\label{VarTFBM3}
{\tilde W}_{\b,c}(t)=\int_{-\infty}^{\infty} \left ((t-u)^{1-\b}\ind{(-c, t)}(u)- ((-u))^{1-\b}\ind{(-c, 0)}(u)\right )^2 (dx).
\end{equation*}
Quite simple calculations give us the following relations
\begin{equation}\label{VarTFBM3A}
{\tilde W}_{\b,c}(t)= \left ((1-\b)^{2} C_0(\b)+(3-2\b)^{-1}\right )t^{3-2\b}+O(t^2),\ \ {\rm as} \ t\to 0,
\end{equation}
\begin{equation}\label{VarTFBM3B}
{\tilde W}_{\b,c}(t)=c^{2(1-\b)}t -\frac{2(1-\b)}{(2-\b)(3-2\b)}c^{3-2\b} \ \ {\rm as}\ t\to \infty.
\end{equation}
In a similar way we can investigate the case $j=9$. We do not provide this case, also the calculations in the case $j=7$ are omitted, since we intend in a separate paper investigate with more details these new  processes, defined in Definition \ref{TFSMIII}.

 Comparing (\ref{VarTFBM3A}), (\ref{VarTFBM3B}) with (\ref{VarZn5}) and (\ref{VarZn6}), we  see that with respect to $t$  they coincide, only there is small difference in constants in these asymptotic.  This difference appears for the reason, that in obtaining $W^{(7)}(t)$ we normalize $S_n^{(7)}(t, X^{(n)})$ in such a way that ${\rm Var}Z_n^{(7)}(t)=1$, while calculating ${\tilde W}_{\b,c}(t)$ we use
$B_{H, c, \s}^{III}(t)$, which is not normalized to have unit variance at point $t=1$ (this can be done by choosing in appropriate way $\s$). Therefore, we can assert that limit processes $U^{(j)}(t), j=7, 9$ are  TFBMIII.



 For a filter $\{a_i\}$ in Theorems \ref{thm1} and \ref{thm3}  we assumed traditional condition for a filter $\b>1/\a$. But since we consider the tapered filter $\{{\bar a}_i^{(n)}\}$ (essentially it is truncated filter), it is possible to consider the case $\b<1/\a$ and even $\b<0$, i.e. the coefficients of a filter unboundedly growing. Since we think that such filters are not realistic, we shall demonstrate only the simple case $\b=0$ in Theorem \ref{thm1}, which has some meaning: filter consists from finite, but big number of coefficients $a_i$ which are  all equal.
 Thus, we consider the linear random process ${\tilde X}^{(n)}=\{ {\tilde X}_{k}^{(n)}=\sum_{j=0}^{\infty} {\tilde a}_{j}^{(n)}\e_{k-j}, \ k\in \bz\}$ with the filter  ${\tilde a}_i^{(n)}=\ind{[0\le i\le \l(n)]}$,  $\l(n)=n^{\g_1}$, and three cases of tapering.
 We prescribe indexes $j=10, 11, 12$ to  cases of strong, weak, and moderate  tapering, respectively. Let us denote ${\tilde Z}_n^{(j)}(t)= \left ( A_n^{(j)}\right )^{-1} S(t,{\tilde X}^{(n)})$.

\begin{prop}\label{prop3}. Let ${\tilde X}_k^{(n)}=\sum_{j=0}^{\infty} {\tilde a}_{j}^{(n)}\e_{k-j}, \ k\in \bz $ be a family of linear processes with the tapered filter ${\tilde a}_i^{(n)}$ and innovations satisfying the condition of Theorem \ref{thm1}. Then, for all $j=10, 11, 12,$, we have
$$
\{Z_n^{(j)}(t)\}\fdd \{U^{(j)}(t)\},
$$
where $U^{(10)}(t)=B(t),  U^{(11)}(t)=B(1)t$, and $ U^{(12)}$ is mean zero Gaussian process with variance, given in (\ref{VarZn11}).
\end{prop}

{\it Sketch of the proof}.
Let  us consider the case $0<\g_1<1$ (strong tapering, the case $j=10$). Repeating the same steps as in the proof of Theorem \ref{thm1} and taking $A_n^{(10)}=n^{\g_1+1/2}$ it is not difficult to get
$$
\lim_{n\to \infty}{\rm Var}Z_n^{(10)}(t)=t .
$$
In the case $\g_1>1$ (weak tapering, the case $j=11$), taking $A_n^{(11)}=n^{1+\g_1/2}$, in a similar way we get
$$
\lim_{n\to \infty}{\rm Var}Z_n^{(11)}(t)=t^2 .
$$
In the case $\g_1=1,  \ \l(n)=cn$ (moderate tapering, the case $j=12$),  we consider the cases $0<t<c$ and $0<c<t$, and the calculations are the same as in the proof of Theorem \ref{thm1}.
 Taking $\left (A_n^{(12)}\right )^2=C_{18}(c)n^{3}$, we get
\begin{equation*}\label{VarZn11}
\lim_{n\to \infty}{\rm Var}Z_n^{(12)}(t)= \left \{\begin{array}{ll}
t^2 C_{19}(t,c) ,               & 0<t\le c, \\
t C_{20}(t,  c) , & \ 0<c<t.
  \end{array} \right.
\end{equation*}
Here $C_{18}(c)=c^2-c^3/3$, for $0<c\le 1$,  $C_{18}(c)=c-1/3$, for $c > 1$,
$C_{19}(t,  c)=c(1-t/(3c))(C_{18}(c))^{-1}$,  for $0<t\le c,$ and for $0<c< t$
$C_{20}(t,  c)=c^2(1-c/(3t))(C_{18}(c))^{-1}.$
Having variances (and at the same covariances) of limit processes for $Z_n^{(j)}(t), j=10, 11, 12,$ it remains to prove the gaussianity of these limit processes. This can be done exactly as in Proposition  \ref{prop2}, by showing $L^{(j)}(2+\d,n, t) \to 0$ as $n\to \infty$ for all $j=10, 11, 12$.
\halmos

 Proposition \ref{prop3} demonstrates that in Theorem \ref{thm1} for tapered filters it is possible to consider the values $\b\le 1/2$. It seems that the following general picture is true. In the case of strong tapering ($0<\g_1<1$) we always get in limit the Brownian motion, despite how big is the growth of  coefficients of a filter (the case of negative $\b$; although in Proposition \ref{prop3} we considered only the case $\b=0,$ but it is easy to see that the same result we get for $\b<0$). This can be explained by the following fact. The variance of $X_k^{(n)}$ grows with  $\b\to -\infty$ and with $\l(n)$, but this growth can be compensated by normalizing. Main factor in proving the asymptotic normality is dependence between summands in a sum $S_n(t, X^{(n)})$, and $X_k^{(n)}, k\ge 1$ for each $n$ are $\l(n)$-dependent. If $\g_1<1$, the number of summands $nt$ grows more rapidly comparing with $\l(n)=n^{\g_1}$, therefore we  get the Brownian motion  as a limit. The same picture can be seen in \cite{Sabzikar}, Theorem 4.3 (i). Contrary, in the case of weak tapering ($\g_1>1$), the main factor becomes dependence between summands, since now $\l(n)$ grows more rapidly comparing with $nt$. In the case of LRD in Theorem \ref{thm1}, the case $j=4$ we have that $H=3/2-\b \to 1$, as $\b\to 1/2$, so one can expect that for all $\b\le 1/2$ the limit process will be degenerate $B(1)t$, as in the case $\b=0$ in Proposition \ref{prop3}. The case $\b<1/2$ is not considered in \cite{Sabzikar}, Theorem 4.3 (ii). More difficult to predict what is happening in the case of moderate tapering,  $\b<1/2$, and our tapering function (\ref{newtaper}).
One may expect that in this case  the limit Gaussian process $U(t)$ will have variance behaving as $t^{2}$ for small $t$ and as $t$ for large $t$, behavior changing at the point $t=c$. Proposition \ref{prop3}, the case $j=12$ with $\b=0$ supports this expectation.

\bigskip

\section{Limit theorems for linear processes with tapered filters and innovations}\label{sec2}

In this section we consider the family of linear processes with tapered filter and innovations, as defined in (\ref{trunkfilinn}). Since for each  tapering there are three cases and there are three cases of dependence, so in total there should be 27 cases. But at  present for linear processes we are able to consider only  hard and soft tapering (leaving intermediate tapering for the future research), therefore it remains 18 cases. Thus,we investigate what happens if in all nine cases, considered in Theorem \ref{thm1} (the Gaussian case) and in Theorem \ref{thm3} (the stable case), we add assumption that innovations are heavy tailed and tapered. As in previous section, we consider the Gaussian and the stable cases separately.

\subsection{The Gaussian case}
 Thus, we consider the family of linear processes with tapered innovations and filter, defined in (\ref{trunkfilinn}). As in Section \ref{sec1}, filter taper is from (\ref{newtaper}), and we assume the same condition (\ref{cond1}) and three types of dependence (SRD, LRD, ND), defined in Section \ref{sec1}. We shall consider the same nine cases  as in Section \ref{sec1}, only now instead of innovations with unit variance we shall consider heavy-tailed tapered innovations. We shall use tapered innovations which  were used in papers \cite{Paul21}, \cite{Paul22}, and \cite{PaulDam5}. Let $\theta =\theta (\a)$ stand for  the standard Pareto distribution with probability density and distribution functions (\ref{paretodf}).
Let  $R$ be the standard exponential random variable with the  density function $e^{-x}$, for $x\ge 0$, independent of $\theta$.
Then  tapered (with tapering parameter $b>1$) standard Pareto random variable $\zeta= \zeta (\a, b)$  can be written as
\begin{equation}\label{tapparetodf}
\zeta (\a, b)=\theta\ind{[\theta<b]}+(b+R)\ind{[\theta\ge b]},
\end{equation}
its density function was written in \cite{Paul22}, see (1.4) therein.
Let us denote $\xi =\xi (\a, b)=\zeta (\a, b)-E\zeta (\a, b)$ . We consider the family (indexed by $n$) of linear random processes ${\bar X}^{(n)}=\{{\bar X}_{k}(b_n),\ k\in \bz \}$,defined in (\ref{trunkfilinn})
where $\xi_{k}(b_n)$ are i.i.d. copies of a random variable $\xi (\a, b_n)$ with a tapering parameter $b_n\to \infty$. (\ref{trunkfilinn}) differs from (\ref{trunkfil})  only in change of innovations. All notations with bar  sign (${\bar Z}_n, {\bar A}_n,$ etc.) mean that ${\bar X}^{(n)}$ is used instead of  $X^{(n)}$.
The reason why we taper standard Pareto random variable instead of more general random variable $\nu$, belonging to the domain of attraction of a stable random variable with exponent $\a$, was explained in Subsection \ref{subsec2}.

Now we recall the notion of the hard and soft tapering, which was introduced in \cite{Chakrabarty} (only with the names of the hard and soft truncation) and was renamed and used in \cite{Paul21}.
\begin{definition}\label{def1} Let  $\{\theta_i \}, \ i\in N$, be i.i.d. random variables with regularly varying tails and with the tail exponent $0<a<2$. A sequence of tapering levels $\{b_n\}$ is called soft, hard, or intermediate tapering for a sequence $\{\theta_i, \ i\in N \}$,  if
$$
\lim_{n\to \infty} n\PP\{|\theta_1|>b_n\}
$$
 is equal to $0$,  $\infty,$ or $0<C<\infty$, respectively.
\end{definition}
 In our case the initial sequence $\{\theta_i\}, \ i\in N,$ is a sequence of standard Pareto random variables, therefore soft and hard tapering is defined as follows: if $b_n=n^\g$ and $\g>1/\a$, we have soft tapering, if $\g=1/\a$ -intermediate tapering, while $0<\g<1/\a$ gives us hard tapering. Clearly, only values $0<\a<2$ are interesting, and in the sequel we shall use this assumption without mentioning.

We consider
$$
{\bar Z}_n(t)={\bar A}_n^{-1} S_n(t, {\bar X}^{(n)}),
$$
where ${\bar X}^{(n)}$ is a family of linear processes with tapered innovations and filter, defined in (\ref{trunkfilinn}). We shall consider the same nine cases, which were considered in Theorem \ref{thm1}, therefore in notations we add index $j$. The following result  shows that in the case of hard tapering of innovations the asymptotic behavior of ${\bar Z}_n^{(j)}(t)$ is the same as of $Z_n^{(j)}(t)$, given in Theorem \ref{thm1}.

\begin{teo}\label{thm2}  For all $j=1,\dots , 9$, we have
$$
\{{\bar Z}_n^{(j)}(t)\}\fdd \{U^{(j)}(t)\},
$$
where   Gaussian limit processes $U^{(j)}(t)$ are defined in Theorem \ref{thm1} and normalizing sequences ${\bar A}_n^{(j)}=A_n^{(j)}(E(\xi_{1}(b(n)))^2)^{1/2}$.
\end{teo}

{\it Proof of Theorem \ref{thm2}}.  As in the proof of Theorem \ref{thm1} there are two steps: calculation of ${\rm Var}S_n(t, {\bar X}^{(n)})$ and proving asymptotic normality of ${\bar Z}_n^{(j)}(t)$.
Since ${\rm Var}S_n(t, {\bar X}^{(n)})={\rm Var}S_n(t, X^{(n)})E(\xi_{1}(b_n))^2$ and ${\bar A}_n^{2}=A_n^{2}E(\xi_{1}(b_n))^2$,  from Proposition \ref{prop1} we get the following result.
\begin{prop}\label{prop4} For all $j=1,\dots , 9$ and $s, t >0$ we have
$$
\lim_{n\to \infty}{\rm Var}{\bar Z}_n^{(j)}(t)= W^{(j)}(t),
$$
$$
\lim_{n\to \infty}{\rm Cov}\left ({\bar Z}_n^{(j)}(t), {\bar Z}_n^{(j)}(s)\right )= W^{(j)}(t)+W^{(j)}(s)-W^{(j)}(|t-s|),
$$
where functions  $W^{(j)}$ are defined in Proposition \ref{prop1} and $({\bar A}_n^{(j)})^2=(A_n^{(j)})^2E(\xi_{1}(b_n))^2$.
\end{prop}
To prove asymptotic normality of ${\bar Z}_n^{(j)}(t)$, we need an analog of Proposition \ref{prop2}. Instead of (\ref{Liapunov}) now we define
$$
{\hat L}^{(j)}(2+\d,n, t):=\frac{\sum_{k=-\infty}^m |d_{n,k,t}^{(j)}|^{2+\d}}{\left (\sum_{k=-\infty}^n (d_{n,k,1}^{(j)})^2\right )^{(2+\d)/2}}\frac{E|\xi_{1}(b_n)|^{2+\d}}{\left (E|\xi_{1}(b_n)|^{2}\right )^{(2+\d)/2}}.
$$
\begin{prop}\label{prop5} For all $j=1, \dots , 9$ there exists $0<\d=\d(j)\le 1$ such that, for $\g<1/\a$, as $n\to \infty$,
\begin{equation}\label{Liapunov3A}
{\hat L}^{(j)}(2+\d,n, t) \to 0.
\end{equation}
In the cases $j=1,2,4,5,7,8$ \ $\d(j)$ can be taken any number in the interval $(0, 1]$, while in the cases $j=3,6,9$ we must take  $0<\d(j)<\min (1, (3-2\b)(\b-1)^{-1})$.
\end{prop}
{\it Proof of Proposition \ref{prop5}}.
Taking into account formula (2.1) in \cite{Paul21} we have
\begin{equation}\label{Liapunov4}
\frac{E|\xi_{1}(b_n)|^{2+\d}}{\left (E|\xi_{1}(b_n)|^{2}\right )^{(2+\d)/2}}\le C(\a, \d)n^{\g \a \d/2}.
\end{equation}
Quantity $\sum_{k=-\infty}^n (d_{n,k,1}^{(j)})^2$ is $(A_n^{(j)})^2$ from Proposition \ref{prop1}, thus it remains to estimate the nominator of the  quantity
$$
{\tilde L}^{(j)}(2+\d,n, t)=\frac{\sum_{k=-\infty}^m |d_{n,k,t}^{(j)}|^{2+\d}}{\left (\sum_{k=-\infty}^n (d_{n,k,1}^{(j)})^2\right )^{(2+\d)/2}}.
$$
In some cases, namely, in the cases $j=2, 4, 5, 7, 8$, we can use the same method of estimation, which was used in Proposition \ref{prop2}:
\begin{equation}\label{Liapunov5}
{\tilde L}^{(j)}(2+\d,n, t)\le \left (\frac{\max_{-\infty<k\le m} |d_{n,k,t}^{(j)}|}{A_n^{(j)}}\right )^{\d}\frac{\sum_{k=-\infty}^m (d_{n,k,t}^{(j)})^2}{\sum_{k=-\infty}^n (d_{n,k,1}^{(j)})^2}.
\end{equation}
Using estimate (\ref{LiapunovA}) and notation (\ref{Liapunov1A}) we see that we need to show that
\begin{equation}\label{Liapunov6}
\left ( \left (A_n^{(j)} \right )^{-1}\max (I_1^{(j)}, I_2^{(j)})n^{\g \a /2}\right)^\d \to 0, \ {\rm as} \ n\to \infty.
\end{equation}
The quantity $\left (A_n^{(j)} \right )^{-1}\max (I_1^{(j)}, I_2^{(j)})$ was estimated in the proof of Proposition \ref{prop2}, therefore it remains to  take these estimates  and compare with $n^{\g \a /2}$. For example, in the case $j=2$ we have $\max (I_1^{(j)}, I_2^{(j)})\le C$ and $A_n^{(2)}$ is of the order $n^{1/2}$, therefore $-1/2 +\g \a /2<0$ if $\g<1/\a$. In a similar way we get (\ref{Liapunov6}) in the cases $j=4, 5, 7, 8$. Thus, we have proved (\ref{Liapunov3A}) for $j=2, 4, 5, 7, 8$, and in these cases $0<\d(j)\le 1$ can be any. But if we try the same method of estimation in the case $j=1,$ we shall get the following result
$$
 \left (A_n^{(j)} \right )^{-1}\max (I_1^{(j)}, I_2^{(j)})n^{\g \a /2}\le C(\a, \b)n^{(-1+2(1-\b)(1-\g_1)+\g\a)/2},
$$
and this quantity tends to zero if $\g\le (1-2(1-\b)(1-\g_1))/\a$. Since in this case $1/2<\b<1$ and $\g_1<1$,  we get more restrictive condition for $\g$. Therefore, we need to estimate ${\tilde L}^{(j)}(2+\d,n, t)$ without using (\ref{Liapunov5}).In the case $j=1$we write
\begin{equation}\label{Jn12}
\sum_{k=-\infty}^m |d_{n,k,t}^{(1)}|^{2+\d}=J_1(n)+J_2(n),
\end{equation}
where (we skip in notations index $j=1$ and use the notations $m=[nt], m_1=[n^{\g_1}]$)
$$
J_1(n)=\sum_{i=0}^\infty \Big | \sum_{k=1}^m a_{k+i}\ind{[0\le k+i\le m_1]} \Big |^{2+\d}, \quad  J_2(n)=\sum_{i=1}^m \Big |\sum_{k=0}^{m-i} a_{k}\ind{[0\le k\le m_1]} \Big |^{2+\d}
$$
In Proposition \ref{prop1} these quantities were estimated in the case $\d=0$, on the other hand, in \cite{Paul22} such sums were estimated with $\d>0$ but with non-tapered filter. Therefore, in the same way we can get
\begin{eqnarray}\label{Jn1}
J_1(n) &=& \sum_{i=0}^\infty \Big | \sum_{k=1}^m a_{k+i}\ind{[0\le k+i\le m_1]} \Big |^{2+\d}=\sum_{i=0}^{m_1} \Big | \sum_{k=1}^{m_1-i} a_{k+i} \Big |^{2+\d} \nonumber \\
       & \le & \int_0^{m_1}\left (\int_0^{m_1-y}(x+y)^{-\b} dx\right )^{2+\d} dy \le C(\b, \d)m_1^{1+(1-\b)(2+\d)},
\end{eqnarray}
\begin{eqnarray}\label{Jn2}
J_2(n) &=& \sum_{i=1}^m \Big | \sum_{k=1}^{\min(m_1, m-i)} a_{k} \Big |^{2+\d}=\sum_{i=1}^{m-m_1} \Big | \sum_{k=1}^{m_1} a_{k} \Big |^{2+\d} +\sum_{i=m-m_1+1}^m \Big | \sum_{k=1}^{ m-i} a_{k} \Big |^{2+\d}  \nonumber \\
       & \le & \int_0^{m-m_1}\left (\int_0^{m_1}x^{-\b} dx\right )^{2+\d} dy + \int_{m-m_1}^m\left (\int_0^{m-y}x^{-\b} dx\right )^{2+\d} dy  \nonumber \\
       & \le & C(\b, \d, t)n^{1+\g_1(1-\b)(2+\d)}.
\end{eqnarray}
From (\ref{Jn12})-(\ref{Jn2}), taking into account that $\g_1<1$ we get
$$
{\hat L}^{(1)}(2+\d,n, t)\le C(\b, \d, t)n^{(-\d+\g\a\d)/2},
$$
and if $\g<1/\a$, then ${\hat L}^{(j)}(2+\d,n, t) \to 0$ as $n\to \infty.$ Thus we have proved (\ref{Liapunov3A}) in the case $j=1$, $\d(1)$ can be chosen arbitrary from the interval $(0, 1)$.
It remains three cases $j=3, 6, 9$, all with ND. Let us consider the case $j=3.$ As in the case $j=1$, we get (again, we skip the index $j=3$)
\begin{equation}\label{Jn31}
J_1(n)   \le m_1^{1+(1-\b)(2+\d)}\int_0^{1}\left (\int_0^{1-y}(x+y)^{-\b} dx\right )^{2+\d} dy,
\end{equation}
\begin{equation}\label{Jn32}
J_2(n)   \le m^{1+(1-\b)(2+\d)}\left ( C(\b, \d)(1+\int_{1-\frac{m_1}{m}}^{1}(1-y)^{(1-\b)(2+\d)}  dy)\right ).
\end{equation}
In the case $j=1$ we had $1/2<\b<1$, now we have $1<\b<3/2$, so we must show that integrals, appearing in (\ref{Jn31}) and (\ref{Jn32}),with appropriate choice of $\d$, are finite. As in \cite{Paul22} (see (2.20)-(2.22) therein) the integral in (\ref{Jn31}) is finite if $(1-\b)(2+\d)>-1$, or $0<\d<(3-2\b)(\b-1)^{-1}$. With the same bounds for $\d$ we get that the integral in (\ref{Jn32}) tends to zero, since $m_1m^{-1} \to 0$ as $n\to \infty.$ Therefore we can take
\begin{equation}\label{delta3}
\d(3)=\min \left (1, \frac{3-2\b}{2(\b-1)}\right )
\end{equation}
and, collecting (\ref{Liapunov4}), (\ref{Jn12}), (\ref{Jn31}), and (\ref{Jn32}), we get for this value of $\d$
$$
{\hat L}^{(3)}(2+\d,n, t)\le C(\b, \d(3), t)n^{(-\d(3)-2(1-\g_1)(\b-1)(2+\d(3)) +\g\a\d(3))/2}.
$$
From this estimate we see that if
$$
\g<\frac{1}{\a}\left (1+\frac{2(1-\g_1)(\b-1)(2+\d(3))}{\d(3)} \right )
$$
then we have (\ref{Liapunov3A}) in the case $j=3$, $\d(3)$ can be taken as in (\ref{delta3}).

Now we consider the case $j=6$, i.e., the case $\g_1>1$ and ND of the filter. Again, skipping the index $j=6$, we can get
\begin{eqnarray*}
J_1(n) &=& \sum_{i=0}^{m_1-m} \Big | \sum_{k=1}^{m} a_{k+i} \Big |^{2+\d} +\sum_{i=m_1-m+1}^{m_1-1} \Big | \sum_{k=0}^{m_1-i} a_{k+i} \Big |^{2+\d} \nonumber \\
       &\le & n^{1+(1-\b)(2+\d)} \int_0^{n^{\g_1-1}-t}\left (\int_0^{t}(x+y)^{-\b} dx\right )^{2+\d} dy \nonumber \\
       &+&  n^{1+(1-\b)(2+\d)} \int_{n^{\g_1-1}-t}^{n^{\g_1-1}} \left (\int_0^{n^{\g_1-1}-y}(x+y)^{-\b} dx\right )^{2+\d} dy   \nonumber \\
\end{eqnarray*}
The integral $ \int_0^{\infty}\left (\int_0^{t}(x+y)^{-\b} dx\right )^{2+\d} dy $ (only with $t=1$) was in the expression of $C_2(\b)$ in the previous section and is finite if $(1-\b)(2+\d)>-1$. It is not very difficult to prove that
$$
\int_{n^{\g_1-1}-t}^{n^{\g_1-1}} \left (\int_0^{n^{\g_1-1}-y}(x+y)^{-\b} dx\right )^{2+\d} dy \to 0 \  \ {\rm as} \ \ n\to \infty.
$$
Therefore, we get
\begin{equation}\label{Jn61}
J_1(n)   \le C(\b, \d, t)n^{1+(1-\b)(2+\d)}.
\end{equation}
Under the same condition $(1-\b)(2+\d)>-1$ we get
\begin{eqnarray*}
J_2(n) &=& \sum_{i=1}^{m} \big |\sum_{k=0}^{m-i} a_{k} \big |^{2+\d}=\sum_{i=1}^{m} \big |\sum_{k=m-i+1}^{\infty} a_{k} \big |^{2+\d} \nonumber \\
     &\le & n^{1+(1-\b)(2+\d)}\int_0^t\left (\int_{t-y}^\infty x^{-\b} dx\right )^{2+\d} dy. \nonumber \\
\end{eqnarray*}
Thus, we have
\begin{equation}\label{Jn62}
J_2(n)   \le C(\b, \d, t)n^{1+(1-\b)(2+\d)}.
\end{equation}
Now taking the value of $\d$ as in (\ref{delta3}) and collecting (\ref{Liapunov4}), (\ref{Jn61}), and (\ref{Jn62}),  we get for this value of $\d$
$$
{\hat L}^{(6)}(2+\d,n, t)\le C(\b, \d, t)n^{(-\d +\g\a\d)/2}.
$$
Thus, we have (\ref{Liapunov3A}) in the case $j=6$, $\d(6)$ can be taken as in (\ref{delta3}). The proof of the last case $j=9$ goes along the same lines as in the case $j=6,$ only now we must consider separately the cases $0<c\le t$ and $0<t<c$ and there will  appear the following integrals:
$$
\int_0^{c-t}\left (\int_0^{t}(x+y)^{-\b} dx\right )^{2+\d} dy, \ \  \int_{c-t}^{c} \left (\int_0^{c-y}(x+y)^{-\b} dx\right )^{2+\d} dy.
$$
Again, for the finiteness of these integrals we need the condition $(1-\b)(2+\d)>-1$, and we get (\ref{Liapunov3A}) in the case $j=9$, $\d(9)$ can be taken as in (\ref{delta3}). Proposition \ref{prop5} is proved.
\halmos
Propositions \ref{prop4} and \ref{prop5} prove Theorem \ref{thm2}.
\halmos

\subsection{The stable case}

We consider the same situation as in the section 2.2, only we taper Pareto innovations, that is, we consider the family of linear processes
\begin{equation}\label{trunkfil2}
{\hat X}^{(n)}=\{{\hat X}_k^{(n)} \}, \quad {\rm where} \ \  {\hat X}_k^{(n)}=\sum_{j=0}^\infty {\bar a}_j^{(n)}\eta_{k-j}(b(n)), \ k\in \bz,
\end{equation}
where tapered filter is defined in (\ref{newtaper}). Tapered innovations are defined in a little bit different way, comparing with  section 3.1.
If  $\theta_b= \theta (\a, b)$ is tapered (with tapering parameter $b>1$) standard Pareto random variable,  defined in (\ref{tapparetodf}), then we define
\begin{equation*}\label{tappar}
{\eta}(b)= {\eta}(\a, b)=\left \{\begin{array}{ll}
              \theta (\a, b), &  {\rm if} \ \ 0<\a<1, \\
              \theta (\a, b)-E\theta (\a, b), & {\rm if}\ \ 1< \a< 2.
               \end{array} \right.
\end{equation*}
Innovations $\eta_i(b_n)$ in (\ref{trunkfil2}) are i.i.d. copies of $\eta(b_n)$ with $b_n=n^\g, \ \g>1/\a$ (hard tapering).
Also we denote $\eta=\theta$, if $0<\a<1$ and $\eta=\theta -E\theta$, if $1<\a<2$. In this section we shall use the easily obtained estimates
\begin{equation}\label{estimate}
E|{\eta}(\a,b)- {\eta}|^r\le\left \{\begin{array}{ll}
             C(\a, r)b^{-(\a-r)}, &  {\rm if} \ \ 0<r<\a<1, \\
              C(\a, r) b^{-(\a-r)},  & {\rm if}\ \ 1<r< \a< 2, \\
              C(\a, r) b^{-(\a-1)r}, & {\rm if}\ \ 0<r<1< \a< 2. \\
               \end{array} \right.
\end{equation}
In Section 2.2 we considered sums $Z_n(t)=A_n^{-1}\sum_{k=0}^{\sv{nt}}X_k^{(n)}$ with linear processes $X_k^{(n)}$ defined in (\ref{linpr2}) (tapered filter and non-tapered Pareto innovations). In Theorem \ref{thm3} we found limit processes for these sums in all nine cases. Now we consider the similar sums, only we take linear process from (\ref{trunkfil2}):
$$
V_n(t)=A_n^{-1}\sum_{k=0}^{\sv{nt}}{\hat X}_k^{(n)}
$$
with the same  normalizing constants $A_n$ as in Theorem \ref{thm3}. We shall show that in all nine cases, considered in Section 2.2, limit processes for $V_n$ coincide with corresponding limit processes for $Z_n$. Namely, we prove the following result.

\begin{teo}\label{thm4} For $V_n(t)$, defined by means of a linear random process (\ref{trunkfil2}) with tapered filter and innovations, for all  $j=1, 2, \dots, 9$,
the following relations hold
\begin{equation*}\label{main1}
\left \{V_{n}^{(j)}(t), \ t \in \bbr_+\right \}\fdd \{U_{j}(t),  \ t \in \bbr_+ \
\end{equation*}
Here  limit processes $U_{j}(t)$ were defined in (\ref{Udef}) and (\ref{Hdef}).
\end{teo}

{\it Proof of Theorem \ref{thm4}}. For the proof we shall use the same method which we were using in \cite{Paul21} and \cite{PaulDam5}. Since the beginning of the proof is the same for all nine cases, we suppress the index $j$ from notations. Clearly, the statement of the theorem will follow if we prove that, for  any $\e>0$ and any fixed $t>0$, as $n \to \infty,$
\begin{equation*}\label{convergP}
\PP\{|V_{n}(t)-Z_{n}(t)|>\e \} \pc 0.
\end{equation*}
This relation will follow if for some $r>0$
$$
E|V_{n}(t)-Z_{n}(t)|^r \to 0.
$$
Writing $Z_n$ and $V_n$ as infinite series of independent random variables, and, as in \cite{PaulDam5}, applying trivial inequality, if $0<r<1$, and, if $1\le r <2$ then applying the well-known result for sums of independent random variables (see 2.6.20 in \cite{Petrov}; we recall that $E(\eta_{i} -\eta_{i}(b_n))=0$ in the case $\a>1$) we get
$$
E|V_{n}(t)-Z_{n}(t)|^r\le C A_n^{-r}\sum_{i=-\infty}^{\sv{nt}}|d_{n, i, t}|^r E|{\eta}_i(\a,b_n)- {\eta}_i|^r.
$$
Thus, we must show  that, for each  $j=1, 2, \dots, 9$, it is possible to choose $0<r=r_j<\a$ such that
\begin{equation}\label{mainest}
I_n^{(j)}:=\left (A_n^{(j)}\right )^{-r_j}\sum_{i=-\infty}^{\sv{nt}}|d_{n, i, t}^{(j)}|^{r_j} E|{\eta}(\a,b_n)- {\eta}|^{r_j} \to 0.
\end{equation}
Since quantities $d_{n, i, t}^{(j)}$ with tapered filter  were investigated and $A_n^{(j)}$ was introduced in Section 2.2, we shall prove (\ref{mainest}) only for several $j$. Let us take $j=1$, this means that we have $1/\a<\b<1, \ \g>1/\a, \ 0<\g_1<1$. Since $1<\a<2$, in (\ref{estimate}) we take $1\le r_1<\a$ and $A_n^{(1)}=n^{1/\a}z_{\b, \g_1, n}^{(1)}=n^{1/\a +\g_1(1-\b)}$. Collecting estimates which we had in Section 2.2 in the case $j=1$, we have
$$
\sum_{i=-\infty}^{\sv{nt}}|d_{n, i, t}^{(1)}|^{r_j}\le C n^{1+\g_1(1-\b)r_1}.
$$
Substituting all these estimates into (\ref{mainest}) with $j=1$ we get
$$
I_n^{(1)}\le Cn^{-(1/\a +\g_1(1-\b))r_1+1+\g_1(1-\b)r_1-\g(\a-r_1) }=Cn^{(\a-r_1)(1/\a-\g)}.
$$
Since $(\a-r_1)(1/\a-\g)<0$ for $\g>1/\a$ and all $1\le r_1<\a$, we have (\ref{mainest}) with $j=1$, and we can take $r_1=1$.

Now  let us take $j=7$ with $\g_1=1, \l (n)=cn$. Again we take $1\le r_7<\a$ and $A_n^{(7)}=n^{1/\a}z_{\b, \g_1, n}^{(7)}=n^{1/\a +(1-\b)}$. Now, as in the Section 2.2, we must consider two cases $t<c$ and $t\ge c$, but since  we do not  need exact asymptotic behavior, it is not difficult to get  estimate
$$
\sum_{i=-\infty}^{\sv{nt}}|d_{n, i, t}^{(7)}|^{r_7}\le C(t, c, \b) n^{1+(1-\b)r_7}.
$$
We see that for $I_n^{(7)}$ we got the same estimate as for $I_n^{(1)}$ with $\g_1=1$, therefore we get (\ref{mainest}) with $j=7$ and we can take $r_7=1$.

The last case with LRD is $j=4$ with $1<\g_1<\infty$. In this case we can proceed as in \cite{PaulDam5}, and, taking $1/\b<r_4<\a$, we should get
$$
I_n^{(4)}\le Cn^{-(1/\a +(1-\b))r_4+1+(1-\b)r_4-\g(\a-r_4) }=Cn^{(\a-r_4)(1/\a-\g)}.
$$
Thus, we have (\ref{mainest}) with $j=4$ and we can take $1/\b<r_4<\a$.

The cases with SRD and ND can be considered in a similar way. For example, in the case $j=2, 0<\g_1<1, \ \g>1/\a, \ \ A:=\sum_{k=0}^\infty |a_k|<\infty$ it is easy to get the estimate
$$
I_n^{(2)}\le Cn^{-r_2/\a +1-\g(\a-r_2) } A^{r_2}=Cn^{(\a-r_2)(1/\a-\g)} A^{r_2}.
$$
We get (\ref{mainest}) with $j=2$ and we can take $r_2=1$. The rest cases we leave  for a reader.
\halmos

\section*{References}

\bibliographystyle{plain}

%


\end{document}